\documentclass[11pt]{article}
\usepackage[a4paper,margin=28mm]{geometry}
\usepackage[T1]{fontenc}
\usepackage[utf8]{inputenc}
\usepackage[english]{babel}
\usepackage[expansion=false]{microtype}
\usepackage{amsmath,amssymb,amsthm,mathtools}
\usepackage[numbers,sort&compress]{natbib}
\usepackage{booktabs}
\usepackage{comment}
\usepackage{tikz,pgfplots}
\usetikzlibrary{arrows.meta}
\pgfplotsset{compat=1.18}
\usetikzlibrary{positioning}
\usepackage[hidelinks]{hyperref}
\numberwithin{equation}{section}
\newtheorem{theorem}{Theorem}[section]
\newtheorem{proposition}[theorem]{Proposition}
\newtheorem{lemma}[theorem]{Lemma}
\newtheorem{corollary}[theorem]{Corollary}
\theoremstyle{remark}
\newtheorem{remark}[theorem]{Remark}
\theoremstyle{definition}
\newtheorem{definition}[theorem]{Definition}
\theoremstyle{plain}
\newcommand{\R}{\mathbb R}
\newcommand{\E}{\mathbb E}
\newcommand{\Prob}{\mathbb P}
\newcommand{\1}{\mathbf 1}
\newcommand{\dd}{\,\mathrm d}
\usepackage{hyperref}
\usepackage{cleveref}

\usepackage{todonotes}

\hypersetup{pdftitle={From random intersection to spatial graphs}}

\title{From random intersection to spatial graphs: a unified asymptotic framework}

\author{
Alessandra Bianchi$^{1}$\thanks{\texttt{alessandra.bianchi@unipd.it}} \qquad
Alberto De Nicolò$^{1}$\thanks{\texttt{alberto.denicolo@studenti.unipd.it}} \qquad
Riccardo Michielan$^{2}$\thanks{\texttt{riccardo.michielan@gssi.it}}\\[0.6em]
{\small $^{1}$Dipartimento di Matematica ``Tullio Levi-Civita'', Università degli Studi di Padova, Padova, Italy}\\[0.2em]
{\small $^{2}$Gran Sasso Science Institute, L'Aquila, Italy}
}

\date{}

\begin{document}
\maketitle

\begin{abstract}
Random intersection graphs and spatial random graphs provide different mechanisms for generating edges and clustering. We investigate their interplay through geometric random intersection graphs, where vertices are adjacent whenever they share a spatially located group. First, a Poisson decomposition by attachment sets yields exact subgraph probabilities and laws of large numbers for connected subgraph densities. It distinguishes triangles supported by a common group from those requiring three distinct groups. Second, explicit representations of global and average local clustering reveal their different dependence on model parameters and the role of shared groups in edge correlations. Finally, by varying the strength and spatial range of group memberships, we establish convergence to random connection, random geometric, random intersection, Erdős–Rényi, and complete graphs. We characterize the corresponding clustering limits and show that clustering need not converge to that of the limiting graph, highlighting the distinction between spatially local convergence and global graph statistics.
\end{abstract}

\noindent\textbf{Keywords:}
geometric random intersection graphs; clustering; subgraph counts; Poisson processes; random geometric graphs.

\section{Introduction}\label{sec:introduction}

Random intersection graphs and spatial random graphs offer two different ways of modelling the formation of edges and triangles. In a random intersection graph, vertices are adjacent when they share an attribute or belong to a common group; each group therefore induces a clique. In spatial random graphs, connections depend on the positions of the vertices, and clustering may arise because two neighbours of a vertex are also likely to be close to one another. These mechanisms are not mutually exclusive. Understanding how they interact, and how one or the other emerges under suitable scalings, is the starting point of this work.

The geometric random intersection graph (GRIG), introduced by Deijfen and Michielan in \cite{deijfen-michielan}, provides a natural setting for this question. Vertices and groups are represented by two independent Poisson point processes in $\mathbb R^d$, with intensities $\lambda$ and $\mu$. Each vertex joins a group independently with a probability determined by their relative positions, and two vertices are connected whenever they share at least one group. For a membership kernel $g$, the probability of an edge between vertices at displacement $x$ is $1-\exp\{-\mu(g*g)(x)\}$. Thus the GRIG has the same pairwise connection probabilities as a suitable random connection model. Nevertheless, its edges are not independent: a single group may generate several edges simultaneously. This distinction becomes essential when studying subgraphs and clustering, and it raises a broader question: to what extent can the same construction reproduce the behaviour of different classes of random graphs?

We address this question from three complementary perspectives: the structure and densities of finite subgraphs, the clustering properties of the infinite graph, and the asymptotic limits of both the graph and its clustering coefficients. The first two describe how geometry and shared groups shape the model at fixed parameters. The third shows how varying the membership mechanism connects it to classical random intersection and spatial graph models, while revealing which structural features are, or are not, preserved by these limits.

To study subgraphs beyond their edge marginals, our first step is an exact decomposition of groups according to their attachments to a prescribed finite set of vertices. By Poisson marking, the numbers of groups attached to each exact subset are independent Poisson random variables. This yields explicit formulas for induced and non-induced subgraph probabilities. Together with Mecke's formula and spatial ergodicity, it also gives almost-sure and $L^1$ laws of large numbers for the subgraph densities in expanding observation windows.

Triangles provide a particularly informative example. We distinguish \emph{clique triangles}, whose three vertices share at least one common group, from \emph{cycle triangles}, which have no common group and require three distinct groups to support their edges. This decomposition allows us to quantify the contribution of each mechanism to triangle formation and identify which one dominates in different regimes.

We then investigate global and average local clustering, deriving conditional representations that reveal their different dependence on vertex intensity and clarify the role of shared groups in clustering and edge correlations. In particular, global clustering turns out to be independent of vertex intensity, whereas average local clustering generally is not.

The main asymptotic analysis concerns the interplay between the strength and spatial range of group memberships. By varying these two scales, we establish convergence in total variation on bounded spatial windows to several classical random graph models: random connection, random geometric, random intersection, Erd\H{o}s--R\'enyi, and complete graphs. A comparison with independent-edge graphs handles several regimes, while the random intersection limit is established by directly coupling group memberships. These results provide a unified framework in which different random graph structures emerge from a common construction under suitable scaling regimes.

Under the same scalings, we characterize the asymptotic behaviour of global and average local clustering. We show that global clustering depends on the scaling parameters through only two quantities, measuring respectively the number of groups joined by a single vertex and shared by two nearby vertices. A central finding is that convergence on bounded spatial windows does not necessarily imply convergence of clustering to that of the limiting graph. While the two limits agree in the fixed-range regimes, they may differ when the interaction range diverges, as illustrated by the Erd\H{o}s--R\'enyi limit.

Finally, we show that even when clustering converges to the value of the random geometric graph, the underlying mechanisms of triangle formation may exhibit different asymptotic behaviours. In particular, in the saturated regime $\mu\to\infty$, the relative contributions of clique and cycle triangles undergo a transition governed by the extent to which groups are shared by three nearby vertices.

\paragraph{Relation to earlier work.}
The GRIG was introduced in \cite{deijfen-michielan}, where degree properties and percolation were investigated.
A closely related hypergraph model was introduced and analysed in \cite{brun26}.
Earlier clustering calculations and simulations on the GRIG appear in \cite{de-nicolo-thesis}.
Subgraph counts and clique-cover representations for non-spatial random intersection graphs were studied in \cite{karonski-scheinerman-singer,rybarczyk-stark}, while related results on clustering and community-based constructions can be found in \cite{bloznelis,deijfen-kets,hofstad-komjathy-vadon}. Clustering in a different spatial, conditionally independent-edge model is considered in \cite{hofstad-hoorn-maitra}. Our Erd\H{o}s--R\'enyi regime is related to equivalence results for non-spatial random intersection graphs \cite{fill-scheinerman-singer,rybarczyk-equivalence}, including quantitative bounds in total variation \cite{kim-lee-na,brennan-bresler-nagaraj}. In our spatial setting, we obtain such bounds for the graph induced by the vertices in any fixed bounded window. Random connection models and geometric graph subgraph counts provide further points of comparison \cite{meester-roy,penrose,bachmann-reitzner,last-nestmann-schulte}. Our proofs use Poisson marking, Mecke's formula, and spatial ergodic theory \cite{last-penrose,NguyenZessin79}. The distinguishing feature of our approach is that it retains the supporting groups throughout the analysis, making it possible to compare not only graph limits but also the mechanisms that generate clustering.

\paragraph{Organization of the paper.}
Section~\ref{sec:model} introduces the GRIG model and 
the Poisson decomposition of group attachments, 
which yields exact formulas for subgraph probabilities 
and the law of large numbers for connected subgraph densities stated in Theorem~\ref{thm:subgraph-density}.
Section~\ref{sec:clustering} studies global and average local clustering, including their clique--cycle decomposition presented in Corollary~\ref{cor:C} and Theorem~\ref{thm:local}, and establishes exact expressions for edge correlations and degree variance. 
It also investigates the dependence of local clustering on vertex intensity, given in Corollary~\ref{cor:local-limits}, and provides the general criterion for clustering convergence 
stated in Theorem~\ref{thm:clustering-criterion}.
Section~\ref{sec:graph-limits} introduces the scaling framework and derives a total variation bound for comparison with independent-edge graphs,
leading to the spatially local convergence results of Theorem~\ref{thm:local-limits}.
Section~\ref{sec:clustering-limits} characterizes global clustering through two parameters, 
establishes its asymptotic limits in Theorem~\ref{thm:boundary}, and compares clustering convergence with the limiting graph models in Corollary~\ref{cor:classical-clustering}. It then explains why spatially local limits and clustering limits may differ, and concludes with the clique--cycle transition under saturation stated in Theorem~~\ref{thm:rgg-clique-cycle}. Section~\ref{sec:conclusion} summarizes the results and discusses open directions.

\section{GRIG model and subgraph counting}\label{sec:model}

In this section, we introduce the model and derive exact formulas for
induced and non-induced subgraph probabilities, by classifying groups
according to the vertices they join. Combined with spatial ergodicity, these
formulas give laws of large numbers for the densities of connected
subgraphs. These results provide the foundation for the clustering analysis and the convergence results developed in the following sections.

\subsection{The model}

Let $V$ and $U$ be independent homogeneous Poisson point processes on $\R^d$, with intensities $\lambda>0$ and $\mu>0$, where $d\ge1$. We call their points vertices and groups, respectively, and identify them with their positions. Let $g:\R^d\to[0,1]$ be measurable, radial and non-increasing, with
\begin{equation}\label{eq:g-assumption}
0<\lVert g\rVert_1:=\int_{\R^d}g(z)\dd z<\infty.
\end{equation}
Conditionally on $V$ and $U$, every vertex--group pair $(v,u)$ is linked independently with probability $g(v-u)$. We write $v\leftrightarrow u$ if vertex $v$ is linked to group $u$. Two vertices are adjacent in the projected graph $G$ if they share at least one group, with the rule that multiple shared groups still produce only one edge. The resulting object is called \emph{Geometric Random Intersection Graph} (GRIG).

Any vertex of $V$ belongs to a Poisson number of groups with mean $\mu\lVert g\rVert_1$, and a group has a Poisson number of members with mean $\lambda\lVert g\rVert_1$. Thus $\mu$ controls how many groups a vertex may join, while $\lambda$ controls how many other vertices can be reached through them.
Write
\begin{equation}\label{eq:f-definition}
f(x):=(g*g)(x)=\int_{\R^d}g(z)g(z-x)\dd z.
\end{equation}
The number of groups shared by two prescribed vertices at displacement $x$ is Poisson with mean $\mu f(x)$ \cite{deijfen-michielan}, and their connection probability is therefore 
\begin{equation}\label{eq:edge-probability}
p(x)=1-e^{-\mu f(x)}.
\end{equation}
In particular, every vertex has finite degree almost surely and the mean degree is $\lambda\int p(z)\dd z\le\lambda\mu\lVert g\rVert_1^2$ (see \cite{deijfen-michielan}). A random connection model with connection function $p$ has the same edge probabilities, but its edges are independent given the vertex positions,
while a random geometric graph is a particular random connection model that connects pairs at a distance less than a fixed radius. The GRIG instead uses the projection rule of a random intersection graph, with spatially varying membership probabilities. 

When studying a typical vertex, we add a vertex at the origin and leave the two Poisson processes otherwise unchanged. This is the Palm construction for the vertex process \cite{last-penrose}. Throughout the paper, expectations and probabilities
 refer to this rooted graph and are denoted respectively by $\E$ and $\mathbb P$, without an additional Palm superscript. 

\subsection{Subgraph probabilities}\label{sec:subgraph_probabilities}

Fix distinct positions $\mathbf v=(v_1,\ldots,v_k)$, where $k\ge2$, 
identify vertex $i$ with $v_i$ and let  $G[\mathbf v]$ be the graph induced on them.
For any $S\subseteq[k]$, we say that a group has type $S$ if it joins exactly the vertices indexed by $S$ among the prescribed tuple, and group memberships outside the tuple are irrelevant. Let $N_S$ count groups of type $S$. We will mainly use types with $|S|\ge2$, since the others do not produce an edge among the selected vertices.

\begin{lemma}[Independent group counts]\label{lem:types}
For fixed distinct positions $\mathbf v$, the random variables $N_S$, $|S|\ge2$, are independent 
and have Poisson distribution with parameter $\Lambda_S(\mathbf v)$, where
\begin{equation}\label{eq:LambdaS}
\Lambda_S(\mathbf v)= \mu \int_{\R^d}\prod_{i\in S}g(z-v_i)\prod_{j\notin S}(1-g(z-v_j))\dd z.
\end{equation}
\end{lemma}
\begin{proof}
 Mark every group by its links to the selected vertices. The probability of mark $S$ at position $z$ is the integrand in \eqref{eq:LambdaS},
 where the first product requires links to $S$, and the second excludes links to the other selected vertices. 
The marking theorem gives independent Poisson processes for distinct marks \cite[Theorem 5.6]{last-penrose}, and integration over $z$ gives their means. 
\end{proof}

For $k=3$, we abbreviate $\Lambda_{\{i,j\}}(\mathbf v)$ and $\Lambda_{\{1,2,3\}}(\mathbf v)$ as $\Lambda_{ij}$ and $\Lambda_{123}$, respectively, whenever the positions $\mathbf v$ are understood. 

A type $S$ group creates all edges $\{v_i,v_j\}$ with $i,j \in S$ in the GRIG. It follows that the set of group types which are present completely determines the induced graph on the selected vertices. Let $\mathcal S(\mathbf v):=\{S\subseteq[k]: |S|\ge2,\ N_S>0\}$ denote the random family of group types present in the configuration. For a given family $\mathcal T$ of types, Lemma~\ref{lem:types} gives
\begin{equation}\label{eq:active-types}
\Prob(\mathcal S(\mathbf v) = \mathcal T)=\prod_{S\in \mathcal T}(1-e^{-\Lambda_S})\prod_{S\notin \mathcal T}e^{-\Lambda_S},
\end{equation}
where products of types range over subsets of $[k]$ with at least two elements. The dependence on $\mathbf v$ is omitted when the positions are fixed.

Next, let $H=([k],E(H))$ be a fixed labelled graph. We distinguish between the event that $H$ occurs as a subgraph of $G[\mathbf v]$ and the event that $G[\mathbf v]$ is exactly $H$. Accordingly, define
\begin{equation}\label{def:prob_subgraph}
p_H(\mathbf v):=\Prob(H\subseteq G[\mathbf v]), \qquad p_H^{\mathrm{ind}}(\mathbf v):=\Prob(H=G[\mathbf v]).
\end{equation}

The decomposition into independent group types yields explicit expressions for both probabilities.

\begin{proposition}[Subgraph probabilities]\label{prop:subgraph-prob}
For any labelled graph $H$ on $[k]$ and all distinct positions $\mathbf v$,
\begin{equation}\label{eq:subgraph-prob}
p_H(\mathbf v)=\sum_{\mathcal T}\Prob(\mathcal S(\mathbf v)=\mathcal T) \1_{\left\{E(H)\: \subseteq\:\bigcup_{S\in\mathcal T}\binom{S}{2}\right\}},
\end{equation}
and
\begin{equation}\label{eq:induced-prob}
p_H^{\mathrm{ind}}(\mathbf v)=\sum_{\mathcal T}\Prob(\mathcal S(\mathbf v)=\mathcal T) \1_{\left\{E(H)\: = \:\bigcup_{S\in\mathcal T}\binom{S}{2}\right\}},
\end{equation}
where $\mathcal T$ ranges over families of subsets of $[k]$ of size at least two, and $\Prob(\mathcal S(\mathbf v)=\mathcal T)$ is explicitly given by \eqref{eq:active-types}.
\end{proposition}

\begin{proof}
Each group of type $S$ creates a clique on the vertices indexed by $S$. Consequently, the edge set of $G[\mathbf v]$ is precisely $\bigcup_{S\in\mathcal S(\mathbf v)}\binom{S}{2}$. The event that $H$ occurs as a subgraph corresponds to this union containing $E(H)$, whereas an induced copy requires equality. Summing over the corresponding disjoint events $\{\mathcal S(\mathbf v)=\mathcal T\}$ gives the two formulas.
\end{proof}

The description of subgraphs through supporting groups is closely related to the clique-cover approach used for random intersection graphs \cite{karonski-scheinerman-singer,rybarczyk-stark,rybarczyk-serafin-subgraph-probabilities}. Indeed, for an induced copy of $H$, the active group types correspond to cliques of $H$ whose edge sets cover $E(H)$. In the spatial Poisson setting, the independent counts of exact attachment types provide explicit, exact expressions for both induced and non-induced subgraph probabilities, while retaining the dependence on vertex positions.

\subsection{Subgraph densities}

Throughout, we use the following spatial averaging principle.
Let $W_n=[-n,n]^d$ and let $\xi(v)$ be a nonnegative measurable
vertex score, covariant under translations of the full bipartite
configuration. We interpret $\E[\xi(0)]$ under the Palm law of
the vertex process, including the group and membership configuration.

The bipartite construction is mixing, hence ergodic, under integer
translations by independence of its induced subgraphs in disjoint
bounded regions and standard local approximation
(cf.\ the proof of \cite[Proposition~8.13]{last-penrose}).
Thus, whenever $\E[\xi(0)]<\infty$, the spatial ergodic theorem
\cite[Corollary~4.20 and Remark~4.28(4)]{NguyenZessin79}
and the Palm--Campbell formula give
\[
\frac{1}{|W_n|}\sum_{v\in V\cap W_n}\xi(v)
\longrightarrow \lambda\E[\xi(0)]
\qquad\text{almost surely and in }L^1.
\]

The same limit holds for rooted counts of a fixed finite subgraph
when only occurrences whose vertices lie in $W_n$ are retained.
Indeed, let $\xi_{W_n}(v)$ denote this restricted score and let
$\xi_L(v)$ count only occurrences within supremum distance
$L\in\mathbb N$ of $v$. For $n>L$,
\[
\sum_{v\in V\cap W_{n-L}}\xi_L(v)
\le
\sum_{v\in V\cap W_n}\xi_{W_n}(v)
\le
\sum_{v\in V\cap W_n}\xi(v).
\]
Dividing by $|W_n|$, applying spatial averaging, and then letting
$L\to\infty$ gives the almost sure limit, since
$\xi_L(0)\uparrow\xi(0)$.
Convergence in $L^1$ follows from uniform integrability,
inherited from the unrestricted averages.

We now apply this principle to subgraph counts in
\[
G|_{W_n}:=G[V\cap W_n],
\]
the induced subgraph on the vertices lying in $W_n$.
Note that only vertices are restricted to the window, while groups outside $W_n$
may still support edges between them.

Let $H$ be a connected graph on $k\ge2$ labelled vertices. For a bounded window $W$, define the number of non-induced copies of $H$ by
\begin{equation}\label{eq:NH}
N_H(W):=\frac{1}{|\operatorname{Aut}(H)|}\sum_{(v_1,\ldots,v_k)\in(V\cap W)^k_{\ne}}\1_{\{H\subseteq G[v_1,\ldots,v_k]\}},
\end{equation}
and define $N_H^{\mathrm{ind}}(W)$ analogously by requiring equality instead of containment. The automorphism factor removes multiple labellings of the same embedding, while distinct embeddings on the same vertex set are still counted.

Recalling the probabilities \(p_H\) and \(p_H^{\mathrm{ind}}\)
from the previous subsection, we define the corresponding rooted
configuration integrals
\begin{equation}\label{eq:IH}
I_H:=\int_{(\R^d)^{k-1}}p_H(0,v_2,\ldots,v_k)\,\dd v_2\cdots\dd v_k,
\end{equation}
and
\begin{equation}\label{eq:IH-induced}
I_H^{\mathrm{ind}}:=\int_{(\R^d)^{k-1}}p_H^{\mathrm{ind}}(0,v_2,\ldots,v_k)\,\dd v_2\cdots\dd v_k,
\end{equation}
where, by translation invariance, the first reference vertex
is placed at the origin, and $v_2,\ldots,v_k$ represent the
positions of the remaining vertices relative to it.
\begin{theorem}[Subgraph densities]\label{thm:subgraph-density}
Under \eqref{eq:g-assumption}, for every connected labelled graph $H$ on $[k]$, $k\ge2$, the configuration integrals $I_H$ and $I_H^{\mathrm{ind}}$ in \eqref{eq:IH} and \eqref{eq:IH-induced} are finite. Moreover,
\begin{equation}\label{eq:subgraph-density}
\frac{N_H(W_n)}{|W_n|}\underset{n\to\infty}{\longrightarrow}\frac{\lambda^k}{|\operatorname{Aut}(H)|}I_H
\end{equation}
almost surely and in $L^1$. The same conclusion holds for induced copies:
\begin{equation}\label{eq:induced-density}
\frac{N_H^{\mathrm{ind}}(W_n)}{|W_n|}\underset{n\to\infty}{\longrightarrow}
\frac{\lambda^k}{|\operatorname{Aut}(H)|}I_H^{\mathrm{ind}}.
\end{equation}
\end{theorem}

\begin{proof}
We first establish the finiteness of $I_H$. 
Let $\mathfrak P(E(H))$ denote the set of partitions of $E(H)$. 
For $\pi=\{B_1,\ldots,B_m\}\in\mathfrak P(E(H))$, let $S_\ell:=\bigcup_{e\in B_\ell}e$, namely the set of vertices incident to at least one edge in \(B_\ell\),
and define
\[
Y_\pi(\mathbf v):=\sum_{(u_1,\ldots,u_m)\in U_{\ne}^{m}}\prod_{\ell=1}^{m}\prod_{i\in S_\ell}\1_{\{v_i\leftrightarrow u_\ell\}}.
\]
Thus, $Y_\pi$ counts assignments of distinct supporting groups to the blocks of $\pi$. Every realization of $H\subseteq G[\mathbf v]$ admits such an assignment for at least one partition, and therefore
\[
p_H(\mathbf v)\le\sum_{\pi\in\mathfrak P(E(H))}\E[Y_\pi(\mathbf v)].
\]

By the multivariate Mecke formula and independence of memberships,
\[
\E[Y_\pi(\mathbf v)]=\mu^m\int_{(\R^d)^m}\prod_{\ell=1}^{m}\prod_{i\in S_\ell}g(u_\ell-v_i)\,\dd u_1\cdots\dd u_m.
\]
Since $H$ is connected, the bipartite incidence graph between the $k$ vertices and the $m$ supporting groups is also connected. Choosing a spanning tree, bounding the factors outside the tree by one, and successively integrating the $k+m-1$ non-root positions gives
\[
\int_{(\R^d)^{k-1}}\E[Y_\pi(0,v_2,\ldots,v_k)]\,\dd v_2\cdots\dd v_k\le\mu^m\|g\|_1^{k+m-1}.
\]
Consequently, by Tonelli's theorem,
\[
I_H\le\sum_{\pi\in\mathfrak P(E(H))}\mu^{|\pi|}\|g\|_1^{k+|\pi|-1}<\infty.
\]
The finiteness of $I_H^{\mathrm{ind}}$ follows from $0\le p_H^{\mathrm{ind}}\le p_H$.

We now apply spatial averaging. Let $\xi_H(v)$ be the rooted count obtained from \eqref{eq:NH} by fixing $v_1=v$ and allowing the remaining vertices to range over all of $V$. This score is nonnegative and translation-covariant, and the multivariate Mecke formula gives
\[
\E[\xi_H(0)]=\frac{\lambda^{k-1}}{|\operatorname{Aut}(H)|}I_H<\infty.
\]
Applying the spatial averaging principle stated above to the
window-restricted counts therefore yields
\[
\frac{N_H(W_n)}{|W_n|}\longrightarrow\lambda\E[\xi_H(0)]=\frac{\lambda^k}{|\operatorname{Aut}(H)|}I_H
\]
almost surely and in $L^1$, proving \eqref{eq:subgraph-density}. The same argument, replacing containment by equality and using $I_H^{\mathrm{ind}}<\infty$, proves \eqref{eq:induced-density}.
\end{proof}

The connectedness assumption is essential for the normalization by spatial volume: disconnected subgraphs may have components that move arbitrarily far apart, producing counts with a different order of growth. For the two subgraphs of primary interest in the next section, namely the wedge $K_{1,2}$ and the triangle $K_3$, we have $|\operatorname{Aut}(K_{1, 2})|=2$ and $|\operatorname{Aut}(K_3)|=6$.

\section{Clustering in the GRIG}\label{sec:clustering}

Now we turn to the clustering properties of the GRIG. Clustering measures the tendency of two neighbours of a vertex to be connected, and is commonly quantified by two related but distinct statistics: the global clustering coefficient and the average local clustering coefficient.

For a finite graph $F$, let $D_F(v)$ denote the degree of a vertex $v$ and $T_F(v)$ the number of triangles containing it. The global clustering coefficient is then given by
\begin{equation}\label{def-Cglob}
C_{\mathrm{glob}}(F):=\frac{\sum_{v\in V(F)}T_F(v)}{\sum_{v\in V(F)}\binom{D_F(v)}{2}},
\end{equation}
while the local clustering coefficient and its average are given by
\begin{equation}\label{def-c-loc}
c_F(v):=\frac{T_F(v)}{\binom{D_F(v)}{2}}, \qquad C_{\mathrm{loc}}(F):=\frac{1}{|V_{\ge2}(F)|}\sum_{v\in V_{\ge2}(F)}c_F(v),
\end{equation}
where $V_{\ge2}(F)$ is the set of vertices in $F$ with degree at least two. We adopt the convention that both coefficients are zero when their respective denominators vanish. 

For the infinite graph $G$, we consider the clustering coefficients on the restrictions $G|_{W_n}$ and their limits as $n\to\infty$.

\subsection{Global clustering}

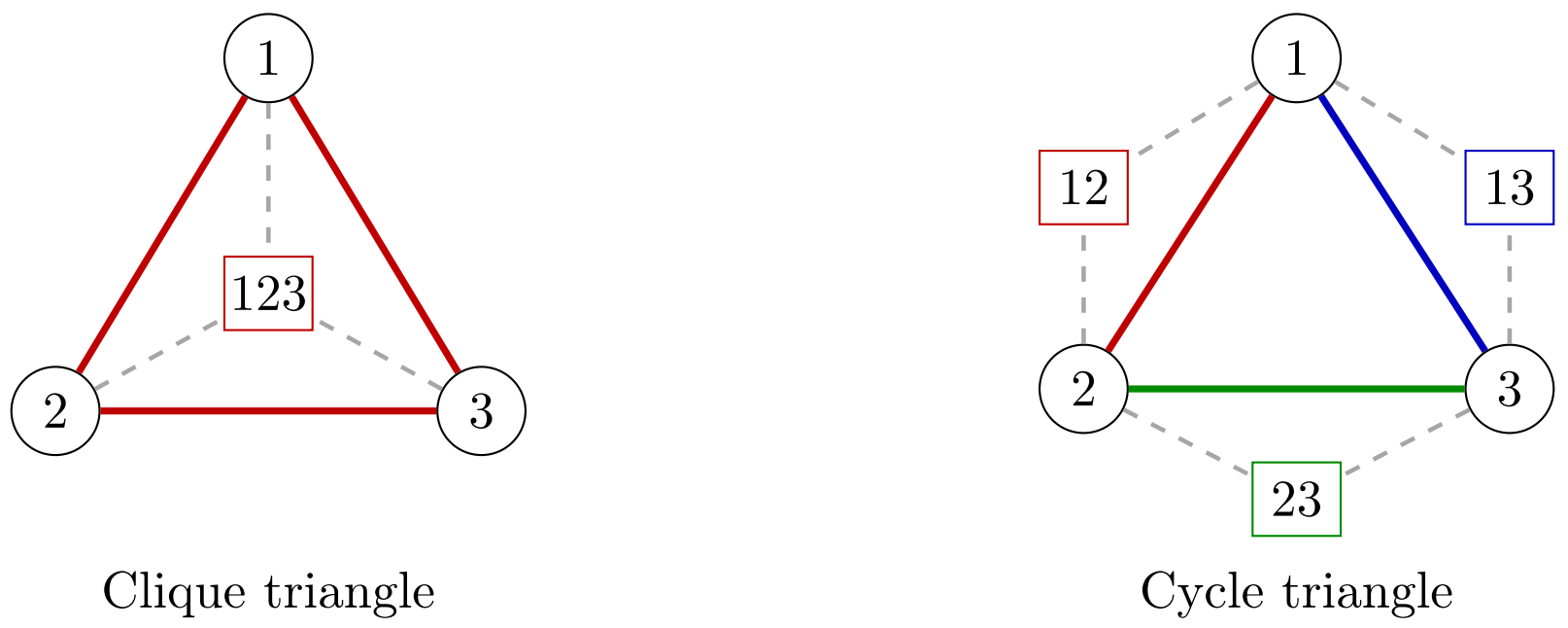
\begin{figure}[tbp]
\centering
\begin{tikzpicture}[
  vertex/.style={circle,draw=black,fill=white,minimum size=6mm,inner sep=1pt},
  group/.style={rectangle,draw,fill=white,minimum height=5mm,minimum width=6mm,inner sep=1pt},
  membership/.style={gray!70,dashed,line width=0.85pt},
  projected/.style={line width=1.35pt},
  every node/.style={font=\small}
]
\begin{scope}[xshift=0cm]
\node[vertex] (a1) at (0,1.6) {$1$};
\node[vertex] (a2) at (-1.45,-.8) {$2$};
\node[vertex] (a3) at (1.45,-.8) {$3$};
\node[group,draw=red!75!black] (a123) at (0,0) {$123$};
\draw[membership] (a1)--(a123) (a2)--(a123) (a3)--(a123);
\draw[projected,red!75!black] (a1)--(a2) (a1)--(a3) (a2)--(a3);
\node at (0,-2.05) {Clique triangle};
\end{scope}
\begin{scope}[xshift=7cm]
\node[vertex] (b1) at (0,1.6) {$1$};
\node[vertex] (b2) at (-1.45,-.65) {$2$};
\node[vertex] (b3) at (1.45,-.65) {$3$};
\node[group,draw=red!75!black] (b12) at (-1.45,.72) {$12$};
\node[group,draw=blue!75!black] (b13) at (1.45,.72) {$13$};
\node[group,draw=green!55!black] (b23) at (0,-1.4) {$23$};
\draw[membership] (b1)--(b12) (b2)--(b12)
                  (b1)--(b13) (b3)--(b13)
                  (b2)--(b23) (b3)--(b23);
\draw[projected,red!75!black] (b1)--(b2);
\draw[projected,blue!75!black] (b1)--(b3);
\draw[projected,green!55!black] (b2)--(b3);
\node at (0,-2.05) {Cycle triangle};
\end{scope}
\end{tikzpicture}
\caption{Schematic bipartite supports and projected edges of a triangle. Circles represent vertices and rectangles represent exact group types relative to the triple. Dashed grey lines indicate vertex--group memberships, while solid coloured lines indicate edges in the projected graph; each colour identifies the group supporting the corresponding edge. In the clique class (left), one common group supports all three edges. In the cycle class (right), three distinct pair groups support one edge each, and no group joins all three vertices. The drawing does not specify spatial positions.}
\label{fig:triangle-supports}
\end{figure}

We begin with the probabilities of wedges and triangles in the GRIG. A \emph{wedge} is a path of length two, consisting of a centre and two distinct neighbours; it is not required to be induced, since the two neighbours may be adjacent. A \emph{triangle} is a triple of vertices with all three edges present. In the GRIG, triangles can arise through two distinct mechanisms. We call a triangle a \emph{clique triangle} if at least one group joins all three vertices. Otherwise, it is a \emph{cycle triangle}, in which case three distinct groups are necessary, each supporting one of the three edges. This classification is exhaustive and disjoint, even when multiple supporting groups coexist; see Figure~\ref{fig:triangle-supports}.

Fix three distinct positions $\mathbf v=(v_1,v_2,v_3)$ and recall the exact group-type means $\Lambda_S(\mathbf v)$ introduced in Section~\ref{sec:model}. The decomposition of group types in Lemma~\ref{lem:types} gives the following probabilities.

\begin{proposition}[Wedge and triangle probabilities]\label{prop:local}
The probability that the three vertices form a wedge centred at $v_2$ is
\begin{equation}\label{eq:wedge-disjoint}
p_\wedge(\mathbf v)=(1-e^{-\Lambda_{123}})+e^{-\Lambda_{123}}(1-e^{-\Lambda_{12}})(1-e^{-\Lambda_{23}}).
\end{equation}
Equivalently,
\begin{equation}\label{eq:wedge-probability}
p_\wedge(\mathbf v)=1-e^{-\mu f(v_1-v_2)}-e^{-\mu f(v_2-v_3)}+e^{-\mu f(v_1-v_2)-\mu f(v_2-v_3)+\Lambda_{123}}.
\end{equation}
The probabilities of a clique triangle and a cycle triangle are
\begin{align}
p_{\mathrm{cl}}(\mathbf v)&=1-e^{-\Lambda_{123}},\label{eq:clique-probability}\\
p_{\mathrm{cyc}}(\mathbf v)&=e^{-\Lambda_{123}}\prod_{1\le i<j\le3}(1-e^{-\Lambda_{ij}}).\label{eq:cycle-probability}
\end{align}
Consequently, the triangle probability is $p_\triangle(\mathbf v)=p_{\mathrm{cl}}(\mathbf v)+p_{\mathrm{cyc}}(\mathbf v)$.
\end{proposition}

\begin{proof}
If a group joins all three vertices, both the wedge and the triangle are present. In the absence of such a group, a wedge centred at $v_2$ requires the exact pair types $\{1,2\}$ and $\{2,3\}$ to be present, whereas a triangle requires all three exact pair types. By Lemma~\ref{lem:types}, the corresponding group counts are independent Poisson random variables, giving \eqref{eq:wedge-disjoint}, \eqref{eq:clique-probability} and \eqref{eq:cycle-probability}. Finally, the identity $\mu f(v_i-v_j)=\Lambda_{ij}+\Lambda_{123}$ yields \eqref{eq:wedge-probability}.
\end{proof}

In the notation of Section~\ref{sec:model}, $p_\wedge=p_{K_{1,2}}$ and $p_\triangle=p_{K_3}$, where $K_{1,2} $ is the path with edges $\{1,2\}$ and $\{2,3\}$.
The clique and cycle decomposition provides additional information about the groups supporting a triangle, which is not captured by the projected graph alone.

We now use the results of \Cref{sec:model} to determine the limiting densities of wedges and triangles. By translation invariance, we fix the centre of the wedge at the origin and write $\mathbf v=(x,0,y)$. We adopt this convention for triples of vertices throughout the remainder of this section. Define
\begin{equation}\label{eq:Iw}
I_\wedge:=\int_{\R^{2d}}p_\wedge(\mathbf v)\dd x\dd y.
\end{equation}
Similarly, for the two triangle mechanisms, let
\begin{equation}\label{eq:Itriangles}
I_{\mathrm{cl}}:=\int_{\R^{2d}}p_{\mathrm{cl}}(\mathbf v)\dd x\dd y,\qquad I_{\mathrm{cyc}}:=\int_{\R^{2d}}p_{\mathrm{cyc}}(\mathbf v)\dd x\dd y,
\end{equation}
and write $I_\triangle := I_{\mathrm{cl}} + I_{\mathrm{cyc}}$.

Let $N_\wedge(W)$ and $N_\triangle(W)$ denote wedge and triangle counts
in a general bounded window $W$, as defined in \eqref{eq:NH}. We also write $N_{\mathrm{cl}}(W)$ and $N_{\mathrm{cyc}}(W)$ for the number of triangles in the two classes, determined by the supporting groups of the full bipartite construction.

\begin{proposition}[Wedge and triangle densities]\label{prop:densities}
Under \eqref{eq:g-assumption}, the configuration integrals $I_\wedge$, $I_{\mathrm{cl}}$ and $I_{\mathrm{cyc}}$ in \eqref{eq:Iw} and \eqref{eq:Itriangles} are finite and $I_\wedge>0$. Moreover,
\begin{align}\label{eq:wedge-density}
\frac{N_\wedge(W_n)}{|W_n|}\underset{n\to\infty}{\longrightarrow}\frac{\lambda^3}{2}I_\wedge,
\qquad
\frac{N_\triangle(W_n)}{|W_n|}\underset{n\to\infty}{\longrightarrow}\frac{\lambda^3}{6}I_\triangle,
\end{align}
almost surely and in $L^1$. The two triangle classes satisfy
\begin{equation}\label{eq:source-densities}
\frac{N_{\mathrm{cl}}(W_n)}{|W_n|}\underset{n\to\infty}{\longrightarrow}\frac{\lambda^3}{6}I_{\mathrm{cl}},\qquad \frac{N_{\mathrm{cyc}}(W_n)}{|W_n|}\underset{n\to\infty}{\longrightarrow}\frac{\lambda^3}{6}I_{\mathrm{cyc}}
\end{equation}
almost surely and in $L^1$.
\end{proposition}

\begin{proof}
The limits \eqref{eq:wedge-density}, and the finiteness of $I_\wedge$ and $I_\triangle$, follow directly from Theorem~\ref{thm:subgraph-density}, applied to $K_{1,2}$ and $K_3$, whose automorphism groups have sizes $2$ and $6$, respectively. Moreover, $I_{\mathrm{cl}}$ and $I_{\mathrm{cyc}}$ are finite since their integrands are bounded by $p_\triangle$.

For the source-separated triangle counts \eqref{eq:source-densities},
the rooted scores counting clique and cycle triangles at a vertex
are nonnegative, translation-covariant, and bounded by the integrable
triangle score. By the Mecke formula, their means are
$\lambda^2I_{\mathrm{cl}}/2$ and $\lambda^2I_{\mathrm{cyc}}/2$.
Applying the spatial averaging principle stated above to the
window-restricted counts and dividing by three yields
\eqref{eq:source-densities}, since each triangle is counted
at its three vertices.

Finally, to establish $I_\wedge>0$, observe that $p_\wedge(\mathbf v)\ge 1-e^{-\Lambda_{123}}$, and Tonelli's theorem gives
\[
\int_{\R^{2d}}\Lambda_{123}\dd x\dd y=\mu\|g\|_1^3>0.
\]
Since $\Lambda_{123}$ is nonnegative and finite, it is positive on a set of positive measure. Hence $I_\wedge>0$.
\end{proof}

Recalling the definition of global clustering, we have
\[
C_{\mathrm{glob}}(G|_{W_n})=\frac{3N_\triangle(W_n)}{N_\wedge(W_n)},
\]
with the ratio set to zero when no wedges are present. Similarly, for $i\in\{\mathrm{cl},\mathrm{cyc}\}$, we define $C_{\mathrm{glob}}^{(i)}(G|_{W_n}):=3N_i(W_n)/N_\wedge(W_n)$, with value zero when the denominator vanishes.
 The preceding density limits immediately give the following result.

\begin{corollary}[Global clustering]\label{cor:C}
Under \eqref{eq:g-assumption}, the global clustering coefficient $C_{\mathrm{glob}}(G|_{W_n})$ and its clique and cycle contributions converge almost surely and in $L^1$ as $n\to\infty$ to
\begin{equation}\label{eq:C}
C_{\mathrm{glob}}:=\frac{I_\triangle}{I_\wedge}=C_{\mathrm{glob}}^{(\mathrm{cl})}+C_{\mathrm{glob}}^{(\mathrm{cyc})},
\qquad C_{\mathrm{glob}}^{(\mathrm{cl})}:=\frac{I_{\mathrm{cl}}}{I_\wedge},\qquad C_{\mathrm{glob}}^{(\mathrm{cyc})}:=\frac{I_{\mathrm{cyc}}}{I_\wedge}.
\end{equation}
The three limits do not depend on the vertex intensity $\lambda$.
\end{corollary}

\begin{proof}
By Proposition~\ref{prop:densities}, the wedge density converges almost surely to a strictly positive constant. Hence the ratio defining global clustering is eventually well-defined, and taking ratios in \eqref{eq:wedge-density}--\eqref{eq:source-densities} gives the asserted limits. The common factor $\lambda^3$ cancels, proving independence from $\lambda$. Since the clustering coefficient and its two contributions take values in $[0,1]$, bounded convergence also gives $L^1$ convergence.
\end{proof}

As a consequence of the above analysis, we also obtain the following natural interpretation of the limiting global clustering coefficient  in terms of a typical vertex.
Let $D(0)$ be the degree of the typical vertex and let $T(0)$ be the number of triangles containing it. By the multivariate Mecke formula,
\[
\E[T(0)]
=
\frac{\lambda^2}{2}I_\triangle,
\qquad
\E\left[\binom{D(0)}{2}\right]
=
\frac{\lambda^2}{2}I_\wedge.
\]
It follows from \eqref{eq:C} that
\begin{equation}\label{eq:C-Palm}
C_{\mathrm{glob}}
=
\frac{\E[T(0)]}
     {\E\left[\binom{D(0)}{2}\right]}.
\end{equation}
This is, in fact, the Palm analogue of the finite-window definition of global clustering \eqref{def-Cglob}.

\subsection{Average local clustering}\label{sec:local}

We now turn to the average local clustering coefficient. Unlike global clustering, which weights vertices proportionally to their number of wedges, average local clustering assigns equal weight to each vertex of degree at least two. Consequently, its limiting Palm representation will involve the expectation of the vertex-wise clustering coefficient, rather than a ratio of expectations as in \eqref{eq:C-Palm}.

For the typical vertex, let $c(0):=T(0)/\binom{D(0)}{2}$ when $D(0)\ge2$, and zero otherwise. We will show that the limiting average local clustering coefficient is $\E[c(0)\mid D(0)\ge2]$ and, unlike global clustering, generally depends on the vertex intensity $\lambda$.

To characterize this quantity, we let $U_0$ denote the set of groups 
linked to the origin and start computing the probability of some basic event,
conditional on $U_0$.

By Poisson thinning, $U_0$ is a Poisson point process with intensity measure $\mu g(u)\dd u$, and hence contains finitely many groups almost surely.
Conditionally on $U_0$, a vertex at position $x$ is a neighbour of the origin if it joins at least one group in $U_0$. Its conditional connection probability is therefore

\[
q(x):=1-\prod_{u\in U_0}(1-g(x-u)).
\]
Since the memberships of distinct vertices are independent, the neighbours of the origin form, conditionally on $U_0$, an inhomogeneous Poisson point process with intensity measure $\lambda q(x)\dd x$. In particular, defining
the $U_0$-measurable random variable
\begin{equation}\label{eq:M}
M:=\int_{\R^d}q(x)\dd x,
\end{equation}
we obtain
\begin{equation}\label{eq:degree-mixture}
D(0)\mid U_0\sim\operatorname{Poisson}(\lambda M).
\end{equation}
Moreover, $M\le\|g\|_1|U_0|<\infty$ almost surely.\\

For $M>0$, let $X_1$ and $X_2$ be two independent samples from the conditional neighbour distribution, whose spatial locations have density $q(x)/M$ and whose memberships to the groups in $U_0$ follow the corresponding conditional distribution. We define the conditional closing probability
\begin{equation}\label{eq:Q}
Q:=\Prob(X_1\sim X_2\mid U_0),
\end{equation}
where adjacency may be supported by groups both inside and outside $U_0$. We set $Q=0$ when $M=0$. Notice that $M$ and $Q$ are determined by $U_0$, and their joint distribution does not depend on the vertex intensity $\lambda$.

To obtain an explicit expression for $Q$, define
\begin{equation}\label{eq:s_e_eta}
s(x,y):=\prod_{u\in U_0}(1-g(x-u)g(y-u))\,,
\quad \eta(x,y):=\mu\int_{\R^d}(1-g(z))g(z-x)g(z-y)\dd z.
\end{equation}
Here $s(x,y)$ is the conditional probability on $U_0$ that two inserted vertices at $x$ and $y$ share no group in $U_0$, while $\eta(x,y)$ is the mean number of groups not linked to the origin that are shared by both vertices.

\begin{lemma}[Clique--cycle decomposition of the closing probability]\label{lem:Q-integral}
For $M>0$,
\begin{equation}\label{eq:Q-integral}
Q=\frac{1}{M^2}\int_{\R^{2d}}\left[1-s(x,y)+\bigl(q(x)q(y)-1+s(x,y)\bigr)(1-e^{-\eta(x,y)})\right]\dd x\dd y.
\end{equation}
Moreover, $Q=Q_{\mathrm{cl}}+Q_{\mathrm{cyc}}$, where
\begin{align*}
Q_{\mathrm{cl}}&:=\frac1{M^2}\int_{\mathbb R^{2d}}\bigl(1-s(x,y)\bigr)\,\dd x\,\dd y,\\ Q_{\mathrm{cyc}}&:=\frac1{M^2}\int_{\mathbb R^{2d}}\bigl(q(x)q(y)-1+s(x,y)\bigr)\bigl(1-e^{-\eta(x,y)}\bigr)\,\dd x\,\dd y
\end{align*}
are the conditional probabilities that $X_1$ and $X_2$ share a group of $U_0$, and that they share none but are adjacent, respectively. We set $Q_{\mathrm{cl}}=Q_{\mathrm{cyc}}=0$ when $M=0$.
\end{lemma}

\begin{proof}
Condition on $U_0$ and fix two locations $x,y\in\R^d$. The membership indicators of the inserted vertices with respect to the groups in $U_0$ are independent Bernoulli variables, with respective probabilities $g(x-u)$ and $g(y-u)$. Thus, the probability that both vertices are neighbours of the origin is $q(x)q(y)$.

The probability that $x$ and $y$ share at least one group in $U_0$ is $1-s(x,y)$. In this case, both vertices are neighbours of the origin and are connected through a common group, giving the clique contribution.

Otherwise, both vertices must be neighbours of the origin without sharing a group in $U_0$, an event of probability $q(x)q(y)-1+s(x,y)$. Conditionally on $U_0$, the groups not linked to the origin form an independent Poisson process with intensity measure $\mu(1-g(z))\dd z$. The number of such groups shared by the two inserted vertices is therefore Poisson with mean $\eta(x,y)$, independently of their memberships to $U_0$. Hence, the probability that at least one such group connects them is $1-e^{-\eta(x,y)}$, giving the cycle contribution.

Combining the two disjoint mechanisms, the joint probability that both inserted vertices are neighbours of the origin and are connected is
\[
1-s(x,y)+\bigl(q(x)q(y)-1+s(x,y)\bigr)(1-e^{-\eta(x,y)}).
\]
Integrating over their locations and dividing by $M^2=\int_{\R^{2d}}q(x)q(y)\dd x\dd y$ yields \eqref{eq:Q-integral}. 
\end{proof}

To account for vertices that support at least one wedge, define
\begin{equation}\label{eq:w-function}
w(t):=1-(1+t)e^{-t},\qquad t\ge0.
\end{equation}
Since $D(0)\mid U_0\sim\operatorname{Poisson}(\lambda M)$, we have $\Prob(D(0)\ge2\mid U_0)=w(\lambda M)$. For $i\in\{\mathrm{cl},\mathrm{cyc}\}$, we denote by $c_i(v)$ and $C_{\mathrm{loc}}^{(i)}(F)$ the corresponding local clustering quantities obtained by counting only triangles of type $i$.

\begin{theorem}[Average local clustering]\label{thm:local}
Under \eqref{eq:g-assumption}, it holds that 
\begin{equation}\label{eq:empirical-local-limit}
\begin{aligned}
&C_{\mathrm{loc}}(G|_{W_n}) \underset{n\to\infty}{\longrightarrow} C_{\mathrm{loc}}(\lambda):=\E[c(0)\mid D(0)\ge2]\\
&C_{\mathrm{loc}}^{(i)}(G|_{W_n})
\underset{n\to\infty}{\longrightarrow} C_{\mathrm{loc}}^{(i)}(\lambda)
:=\E[c_ i (0)\mid D(0)\ge2],
\qquad  i \in\{\mathrm{cl},\mathrm{cyc}\},
\end{aligned}
\end{equation}
almost surely and in $L^1$, with
$C_{\mathrm{loc}}(\lambda)
=C_{\mathrm{loc}}^{(\mathrm{cl})}(\lambda)
+C_{\mathrm{loc}}^{(\mathrm{cyc})}(\lambda) \in [0,1]$.
Moreover, the limit and its two contributions satisfy
\begin{equation}\label{eq:local-global}
C_{\mathrm{loc}}(\lambda)
=\frac{\E[Qw(\lambda M)]}
        {\E[w(\lambda M)]},\qquad
C_{\mathrm{loc}}^{(i)}(\lambda)
=\frac{\E[Q_ i  w(\lambda M)]}
        {\E[w(\lambda M)]},
\qquad   i \in\{\mathrm{cl},\mathrm{cyc}\}.
\end{equation}
\end{theorem}

\begin{proof}
We first establish \eqref{eq:empirical-local-limit}.
Applying the spatial averaging principle stated above to the bounded
translation-covariant scores $c(v)$ and
$\1_{\{D(v)\ge2\}}$, we obtain
\[
\frac{1}{|W_n|}
\sum_{v\in V\cap W_n}c(v)
\longrightarrow\lambda\E[c(0)]
\]
and
\[
\frac{1}{|W_n|}
\sum_{v\in V\cap W_n}\1_{\{D(v)\ge2\}}
\longrightarrow\lambda\Prob(D(0)\ge2)
\]
almost surely and in $L^1$.

We next compare the full-graph scores with their finite-window
counterparts. For $v\in V\cap W_n$, let $c_{G|_{W_n}}(v)$ and
$D_{G|_{W_n}}(v)$ denote its local clustering coefficient and degree
in the restricted graph, respectively.
Fix an integer $L\ge1$ and let $A_L(v)$
be the event that $v$ has a neighbour at
$\ell^\infty$-distance greater than $L$.
If $v\in W_{n-L}$ and $A_L(v)$ does not occur,
all its neighbours belong to $W_n$, so its degree
and local clustering coefficient are unchanged
by restriction.
Since both $c(v)$ and
$c_{G|_{W_n}}(v)$ belong to $[0,1]$,
\[
\begin{aligned}
&\frac{1}{|W_n|}
\sum_{v\in V\cap W_n}
|c_{G|_{W_n}}(v)-c(v)|
\le
\frac{|V\cap(W_n\setminus W_{n-L})|}{|W_n|}
+
\frac{1}{|W_n|}
\sum_{v\in V\cap W_n}\1_{A_L(v)}.
\end{aligned}
\]
For fixed $L$, the first term tends to zero almost surely
and in $L^1$, while the spatial averaging principle stated above,
applied to the bounded translation-covariant score
$\1_{A_L(v)}$, shows that the second term converges to
$\lambda\Prob(A_L(0))$ almost surely and in $L^1$.
Since $D(0)<\infty$ almost surely,
$\Prob(A_L(0))\to0$ as $L\to\infty$.
Taking the limit superior as $n\to\infty$ for each fixed
integer $L$, and then letting $L\to\infty$, proves that
the normalized discrepancy vanishes almost surely.
The same conclusion in $L^1$ follows by taking
expectations in the preceding bound.
The same argument applies to the indicators of the
events that a vertex has degree at least two.
Consequently,
\[
\frac{1}{|W_n|}
\sum_{v\in V\cap W_n}c_{G|_{W_n}}(v)
\longrightarrow\lambda\E[c(0)]
\]
and
\[
\frac{|\{v\in V\cap W_n:D_{G|_{W_n}}(v)\ge2\}|}
     {|W_n|}
\longrightarrow\lambda\Prob(D(0)\ge2)
\]
almost surely and in $L^1$.
The latter limit is strictly positive. Indeed,
$\Prob(U_0\ne\varnothing)>0$, and
$M\ge\|g\|_1>0$ on this event. Hence
\[
\Prob(D(0)\ge2)=\E[w(\lambda M)]>0.
\]
Taking the ratio, therefore, gives
\[
C_{\mathrm{loc}}(G|_{W_n})
\longrightarrow
\frac{\E[c(0)]}{\Prob(D(0)\ge2)}
=
\E[c(0)\mid D(0)\ge2]
=
C_{\mathrm{loc}}(\lambda),
\]
where we used $c(0)=0$ when $D(0)<2$.
Since $0\le C_{\mathrm{loc}}(G|_{W_n})\le1$,
bounded convergence also yields convergence in $L^1$.

The same averaging and boundary estimates apply to
$c_ i (v)$, since $0\le c_ i (v)\le c(v)\le1$.
Indeed, whenever all neighbours of $v$ belong to $W_n$,
its incident triangles and their source classes are
unchanged by restriction, since groups remain unrestricted.
This proves the two source-specific limits in
\eqref{eq:empirical-local-limit}, almost surely and in $L^1$.
Moreover, $c=c_{\mathrm{cl}}+c_{\mathrm{cyc}}$ gives
\[
C_{\mathrm{loc}}(\lambda)
=
C_{\mathrm{loc}}^{(\mathrm{cl})}(\lambda)
+
C_{\mathrm{loc}}^{(\mathrm{cyc})}(\lambda)\in[0,1].
\]

We now prove \eqref{eq:local-global}.
Condition on $U_0$ and on $D(0)=k\ge2$.
By the conditional Poisson construction, the $k$ neighbours
have independent locations and membership marks to $U_0$
sampled from the normalized neighbour distribution,
and each unordered pair is connected with probability $Q$.
Consequently,
\[
\E[c(0)\mid U_0,D(0)=k]=Q.
\]
Since $c(0)=0$ when $D(0)<2$, it follows that
\[
\E[c(0)\mid U_0]=Qw(\lambda M).
\]
Taking expectations and using
$\Prob(D(0)\ge2)=\E[w(\lambda M)]$ gives
\[
C_{\mathrm{loc}}(\lambda)
=
\frac{\E[Qw(\lambda M)]}{\E[w(\lambda M)]}.
\]
The same conditional calculation, with $c$ and $Q$
replaced by $c_ i $ and $Q_ i $, respectively,
gives the remaining identities in \eqref{eq:local-global}.
\end{proof}

The same conditional description also gives a representation
of global clustering:
\begin{equation}\label{eq:global-group-representation}
C_{\mathrm{glob}}
=\frac{\E[QM^2]}{\E[M^2]},
\qquad
C_{\mathrm{glob}}^{(i)}=\frac{\E[Q_ i  M^2]}{\E[M^2]},
\quad  i \in\{\mathrm{cl},\mathrm{cyc}\}.
\end{equation}
Indeed, the conditional Poisson degree distribution and the
closing probabilities yield
\[
\E\left[\left.\binom{D(0)}{2}\right|U_0\right]
=\frac{\lambda^2}{2}M^2,
\qquad
\E[T(0)\mid U_0]
=\frac{\lambda^2}{2}QM^2.
\]
Substitution into \eqref{eq:C-Palm} proves the first identity;
the same argument with the source-specific triangle counts
gives the second. The expectations are finite because
$M\le\|g\|_1|U_0|$ and $|U_0|$ is Poisson,
and the denominator is strictly positive since
$\Prob(M>0)>0$.

The representations in \eqref{eq:local-global} and \eqref{eq:global-group-representation} show that both clustering coefficients average the same conditional closing probability $Q$, but with different weights. Average local clustering weights each group environment by the probability $w(\lambda M)$ that the typical vertex has at least two neighbours, whereas global clustering weights it proportionally to $M^2$, its expected number of neighbour pairs. In particular, only the former depends on the vertex intensity $\lambda$. 
 \begin{remark}[Clustering of vertices of given degree] \label{rem:degree-clustering}
The proof of Theorem~\ref{thm:local} also gives the average clustering of the vertices of given degree, often called the clustering spectrum. 
Explicitly, for $k\ge 2$, this is given by 
\begin{equation}\label{eq:degree-clustering}
    \bar c_\lambda(k):=\E[c(0)\mid D(0)=k]=\frac{\E[QM^ke^{-\lambda M}]}{\E[M^ke^{-\lambda M}]}.
\end{equation}

Indeed, given $U_0$ and $D(0)=k\ge2$, the neighbours of the origin are independent with density $q/M$, so that $\E[c(0)\mid U_0,D(0)=k]=Q$,  while
$\Prob(D(0)=k\mid U_0)=e^{-\lambda M}(\lambda M)^k/k!$ by \eqref{eq:degree-mixture}. Both clustering coefficients are averages of $\bar c_\lambda(k)$,
and  writing $p_k:=\Prob(D(0)=k)$, we obtain

 \begin{equation*}
    C_{\mathrm{loc}}(\lambda)=
    \frac{\sum_{k\ge2}p_k\bar c_\lambda(k)}
        {\sum_{k\ge2}p_k},\qquad
    C_{\mathrm{glob}}=\frac{\sum_{k\ge2}k(k-1)p_k\bar c_\lambda(k)}
        {\sum_{k\ge2}k(k-1)p_k}.
\end{equation*}
Thus global clustering gives greater weight to vertices of higher degree; in particular, $C_{\mathrm{glob}}\le C_{\mathrm{loc}}(\lambda)$ whenever
 $\bar c_\lambda$ is non-increasing.
\end{remark}

The group structure underlying clustering also induces dependencies between incident edges. In particular, groups shared by three vertices contribute both to clique triangles and to positive edge correlations. The following proposition makes this connection explicit and shows how these correlations affect the degree variance.

\begin{proposition}[Edge correlations and degree variance]\label{prop:edge-correlations}
For three distinct positions $\mathbf v=(x,0,y)$, we have
\begin{equation}\label{eq:non-edge-correlation}
\Prob(0\nsim x,\ 0\nsim y)=\Prob(0\nsim x)\,\Prob(0\nsim y)\,e^{\Lambda_{123}},
\end{equation}
and consequently
\[
\operatorname{Cov}\bigl(\1_{\{0\sim x\}},\1_{\{0\sim y\}}\bigr)=e^{-\mu f(x)-\mu f(y)}\bigl(e^{\Lambda_{123}}-1\bigr)\ge0.
\]
Moreover, the degree of the typical vertex has mean and variance
\begin{align}\label{eq:overdispersion-E}
\E[D(0)]&=\lambda\int_{\R^d}\bigl(1-e^{-\mu f(x)}\bigr)\dd x,\\
\operatorname{Var}(D(0))&=\E[D(0)]+\lambda^2\int_{\R^{2d}}e^{-\mu f(x)-\mu f(y)}\bigl(e^{\Lambda_{123}}-1\bigr)\dd x\dd y,\label{eq:overdispersion-V}
\end{align}
where the double integral equals $\operatorname{Var}(M)$.
\end{proposition}

\begin{proof}
The events $\{0\nsim x\}$ and $\{0\nsim y\}$ occur together if and only if no group of exact type $\{1,2\}$, $\{2,3\}$ or $\{1,2,3\}$ is present, which by Lemma~\ref{lem:types} has probability $e^{-\Lambda_{12}-\Lambda_{23}-\Lambda_{123}}$. Since $\Lambda_{12}+\Lambda_{123}=\mu f(x)$ and $\Lambda_{23}+\Lambda_{123}=\mu f(y)$, this probability equals $e^{-\mu f(x)-\mu f(y)+\Lambda_{123}}$, while $\Prob(0\nsim x)=e^{-\mu f(x)}$. This proves \eqref{eq:non-edge-correlation}, and the covariance of two indicators equals that of their complements.

For \eqref{eq:overdispersion-E} and \eqref{eq:overdispersion-V}, the conditional Poisson law \eqref{eq:degree-mixture} and the law of total variance give
\[
\operatorname{Var}(D(0))=\E[\lambda M]+\operatorname{Var}(\lambda M)=\E[D(0)]+\lambda^2\operatorname{Var}(M).
\]
Conditionally on $U_0$, vertices at $x$ and $y$ are neighbours of the origin independently, with probabilities $q(x)$ and $q(y)$. Hence $\E[q(x)q(y)]=\Prob(0\sim x,\,0\sim y)$ and $\E[q(x)]=\Prob(0\sim x)$. Since $\E[M^2]<\infty$, Tonelli's theorem gives $\operatorname{Var}(M)=\int\!\!\int\operatorname{Cov}(q(x),q(y))\dd x\dd y$, and the first part gives the integral.
\end{proof}

For a random connection model, the degree of a typical vertex is Poisson distributed, so that its variance equals its mean. In the GRIG, the index of dispersion is
\[
\frac{\operatorname{Var}(D(0))}{\E[D(0)]}=1+\lambda\,\frac{\int_{\R^{2d}}e^{-\mu f(x)-\mu f(y)}\bigl(e^{\Lambda_{123}}-1\bigr)\dd x\dd y}{\int_{\R^d}p(x)\dd x}.
\]
Thus, for fixed $\mu$ and $g$, the excess dispersion is linear in the vertex intensity $\lambda$ and is determined by the integrated covariance between incident edges. This positive dependence and the presence of clique triangles originate from the same mechanism: groups joining three vertices, whose mean number $\Lambda_{123}$ also determines the clique probability $p_{\mathrm{cl}}=1-e^{-\Lambda_{123}}$.

\subsection{Dependence on vertex intensity and asymptotic clustering}

As $\lambda\downarrow0$ the weight $w(\lambda M)$ becomes proportional to $M^2$, as in global clustering, while as $\lambda\to\infty$ it tends to $\1_{\{M>0\}}$. This leads to the following limits.

\begin{corollary}[Sparse- and dense-vertex limits]\label{cor:local-limits}
Under \eqref{eq:g-assumption}, for fixed $\mu$ and $g$, as functions of the vertex intensity,
\begin{equation}\label{eq:local-limits}
\lim_{\lambda\downarrow0}C_{\mathrm{loc}}(\lambda)=C_{\mathrm{glob}},\qquad \lim_{\lambda\to\infty}C_{\mathrm{loc}}(\lambda)=\E[Q\mid M>0].
\end{equation}
\end{corollary}

\begin{proof}
By definition \eqref{eq:w-function} of $w(t)$, it holds that
\[
0\le\frac{w(\lambda M)}{\lambda^2}\le\frac{M^2}{2},
\qquad
\lim_{\lambda\downarrow0}\frac{w(\lambda M)}{\lambda^2}=\frac{M^2}{2} .
\]
Since $\E[M^2]<\infty$ and $0\le Q\le1$, dominated convergence, applied separately to the numerator and denominator of the first identity of \eqref{eq:local-global},
 gives
\[
\lim_{\lambda\downarrow0}C_{\mathrm{loc}}(\lambda)=\frac{\E[QM^2]}{\E[M^2]}=C_{\mathrm{glob}}.
\]
On the other hand, 
\[
\lim_{\lambda\to\infty}w(\lambda M)=\1_{\{M>0\}}.
\]
Since $0\le w\le1$ and $\Prob(M>0)>0$, dominated convergence yields
\[
\lim_{\lambda\to\infty}C_{\mathrm{loc}}(\lambda)=\frac{\E[Q\1_{\{M>0\}}]}{\Prob(M>0)}=\E[Q\mid M>0].
\]
\end{proof}

Note that the difference between the two limiting values admits the covariance representation
\begin{equation}\label{eq:covariance}
\lim_{\lambda\downarrow0}C_{\mathrm{loc}}(\lambda)-\lim_{\lambda\to\infty}C_{\mathrm{loc}}(\lambda)=\frac{\operatorname{Cov}(Q,M^2\mid M>0)}{\E[M^2\mid M>0]}.
\end{equation}
Thus, the sparse- and dense-vertex limits coincide precisely when $Q$ and $M^2$ have zero conditional covariance given $M>0$. In particular, a nonzero covariance implies that $C_{\mathrm{loc}}(\lambda)$ is not constant in $\lambda$.

The preceding result describes the dependence of average local clustering on the vertex intensity $\lambda$, while keeping $\mu$ and $g$ fixed. We now consider the more general setting in which these parameters may vary simultaneously. The following convergence criterion, based on the conditional representations \eqref{eq:local-global} and \eqref{eq:global-group-representation}, will be useful in the analysis of the asymptotic regimes in Section~\ref{sec:clustering-limits}. As the GRIG parameter $\mu$ varies, the typical neighbourhood mass may converge to a positive finite limit, vanish, or diverge. We therefore introduce a deterministic normalization $b_\mu>0$ and investigate the joint convergence of $(M_\mu/ b_\mu,Q_\mu)$.
Here and below, we use the subscript $\mu$ to indicate dependence on this parameter whenever needed. The following criterion covers all three possible limiting regimes of the effective neighbourhood scale $\lambda_\mu b_\mu$, where we allow the vertex intensity $\lambda_\mu$ to depend on $\mu$.

\begin{theorem}[General criterion for clustering convergence]\label{thm:clustering-criterion}
Consider a family of GRIGs with vertex intensities $\lambda_\mu>0$, and let $M_\mu,Q_\mu$ denote the corresponding conditional neighbourhood mass and closing probability. Suppose that, along a limiting regime of $\mu$, there exists a deterministic normalization $b_\mu>0$ such that
\begin{equation}\label{eq:clustering-joint-convergence}
\left(\frac{M_\mu}{b_\mu},Q_\mu\right)\overset{d}{\longrightarrow}(M_*,Q_*),
\end{equation}
where $M_*$ is finite almost surely, $\Prob(M_*>0)>0$, and $0\le Q_*\le1$.

If the family $\{(M_\mu/b_\mu)^2\}$ is uniformly integrable, then
\begin{equation}\label{eq:clustering-global-limit}
C_{\mathrm{glob}}(\mu)\longrightarrow
\frac{\E[Q_*M_*^2]}{\E[M_*^2]}.
\end{equation}

For the average local clustering coefficient, the following conclusions hold.

\begin{enumerate}
\item If $\lambda_\mu b_\mu\to0$ and $\{(M_\mu/b_\mu)^2\}$ is uniformly integrable, then
\begin{equation}\label{eq:clustering-local-small}
C_{\mathrm{loc}}(\lambda_\mu,\mu)\longrightarrow
\frac{\E[Q_*M_*^2]}{\E[M_*^2]}.
\end{equation}

\item If $\lambda_\mu b_\mu\to\rho\in(0,\infty)$, then
\begin{equation}\label{eq:clustering-local-intermediate}
C_{\mathrm{loc}}(\lambda_\mu,\mu)\longrightarrow
\frac{\E[Q_*w(\rho M_*)]}{\E[w(\rho M_*)]}.
\end{equation}

\item If $\lambda_\mu b_\mu\to\infty$ and
\begin{equation}\label{eq:clustering-small-mass-condition}
\lim_{\varepsilon\downarrow0}\limsup_{\mu}
\E\left[
w(\lambda_\mu M_\mu)
\1_{\{M_\mu/b_\mu\le\varepsilon\}}
\right]=0,
\end{equation}
then
\begin{equation}\label{eq:clustering-local-large}
C_{\mathrm{loc}}(\lambda_\mu,\mu)\longrightarrow
\E[Q_*\mid M_*>0].
\end{equation}
Condition \eqref{eq:clustering-small-mass-condition} is automatically satisfied if $\Prob(M_*>0)=1$.
\end{enumerate}
\end{theorem}

\begin{proof}
Write $X_\mu=M_\mu/b_\mu$ so that, by \eqref{eq:global-group-representation},
\[
C_{\mathrm{glob}}(\mu)
=\frac{\E[Q_\mu X_\mu^2]}{\E[X_\mu^2]}.
\]
The joint convergence in distribution \eqref{eq:clustering-joint-convergence}
and uniform integrability of $X_\mu^2$ imply
\[
\E[X_\mu^2]\longrightarrow\E[M_*^2],
\qquad
\E[Q_\mu X_\mu^2]\longrightarrow\E[Q_*M_*^2].
\]
Since $\Prob(M_*>0)>0$, the limiting denominator is strictly positive, proving \eqref{eq:clustering-global-limit}.

For the analysis of the average local clustering, by \eqref{eq:local-global} 
we write
\[
C_{\mathrm{loc}}(\lambda_\mu,\mu)
=\frac{\E[Q_\mu w(\lambda_\mu b_\mu X_\mu)]}
{\E[w(\lambda_\mu b_\mu X_\mu)]}.
\]
Suppose first that $\lambda_\mu b_\mu\to0$. Since $w(x)\sim x^2/2$ as $x\downarrow0$, and $0\le2w(x)\le x^2$ for every $x\ge0$, joint convergence and uniform integrability yield
\[
\frac{2\E[w(\lambda_\mu b_\mu X_\mu)]}{(\lambda_\mu b_\mu)^2}
\longrightarrow\E[M_*^2],
\qquad
\frac{2\E[Q_\mu w(\lambda_\mu b_\mu X_\mu)]}{(\lambda_\mu b_\mu)^2}
\longrightarrow\E[Q_*M_*^2],
\]
and taking the ratio proves \eqref{eq:clustering-local-small}.

If $\lambda_\mu b_\mu\to\rho\in(0,\infty)$, then $(\lambda_\mu b_\mu X_\mu,Q_\mu)\to(\rho M_*,Q_*)$ in distribution. Since $w$ is bounded and continuous, this gives
\[
\E[w(\lambda_\mu b_\mu X_\mu)]\longrightarrow\E[w(\rho M_*)],
\]
and the analogous convergence for the numerator. The limiting denominator is positive, proving \eqref{eq:clustering-local-intermediate}.

Finally, suppose that $\lambda_\mu b_\mu\to\infty$. On every interval $[\varepsilon,\infty)$ with $\varepsilon>0$, $w(\lambda_\mu b_\mu x)\to1$ uniformly. Joint convergence therefore gives convergence of the numerator and denominator after excluding the region $X_\mu\le\varepsilon$, and condition \eqref{eq:clustering-small-mass-condition} ensures that the excluded contribution vanishes as $\varepsilon\downarrow0$. Consequently,
\[
\E[w(\lambda_\mu M_\mu)]\longrightarrow\Prob(M_*>0),
\qquad
\E[Q_\mu w(\lambda_\mu M_\mu)]
\longrightarrow\E[Q_*\1_{\{M_*>0\}}],
\]
which proves \eqref{eq:clustering-local-large}. If $M_*>0$ almost surely, the small-mass condition follows directly from joint convergence, since the probability that $X_\mu$ lies arbitrarily close to zero vanishes.
\end{proof}

When the conditional environment concentrates, the theorem takes a simple form.

\begin{corollary}[Concentrated environments]\label{cor:concentrated}
    In the setting of Theorem~\ref{thm:clustering-criterion}, suppose that $M_\mu/b_\mu\to m$ in $L^2$ and $Q_\mu\to q$ in probability, for constants $m>0$ and $q\in[0,1]$. Then $C_{\mathrm{glob}}(\mu)\to q$ and $C_{\mathrm{loc}}(\lambda_\mu,\mu)\to q$, for every choice of the vertex intensities $\lambda_\mu$.
\end{corollary}

\begin{proof}
    Convergence in $L^2$ implies uniform integrability of $(M_\mu/b_\mu)^2$, and the joint convergence \eqref{eq:clustering-joint-convergence} holds with the constants $M_*=m$ and $Q_*=q$. Hence $C_{\mathrm{glob}}(\mu)\to qm^2/m^2=q$. Every subsequence has a further subsequence along which $\lambda_\mu b_\mu$ converges in $[0,\infty]$. Along it, each of the three cases of Theorem~\ref{thm:clustering-criterion} gives the limit $q$, because $Q_*=q$ is deterministic and $M_*=m>0$, so that \eqref{eq:clustering-small-mass-condition} holds automatically. Therefore $C_{\mathrm{loc}}(\lambda_\mu,\mu)\to q$ along the whole family.
\end{proof}

\section{Spatially local convergence}\label{sec:graph-limits}

We now investigate how the GRIG approaches several limiting random graph models under suitable parameter scalings. The GRIG combines two sources of dependence: spatial geometry, which determines, through locations, the probabilities of vertex-group memberships, and the mediation of edges through shared groups, which introduces dependencies between edges. We show
that, by weakening or removing these mechanisms, one obtains random intersection graphs (RIG), random connection models (RCM), including random geometric graphs (RGG) as a particular case, Erd\H{o}s--R\'enyi graphs (ER), and the complete graph.

Our results concern convergence in total variation of the graph restricted to any fixed, bounded spatial window, while retaining both vertex positions and their connections. We first compare the GRIG with a graph with independent edges, using the group-type decomposition of Section~\ref{sec:model}, and then apply this comparison to all limiting models except the random intersection graph, which is obtained by a direct coupling of memberships. Throughout this section, the vertex intensity $\lambda$ is kept fixed, while the group intensity $\mu$ and the membership kernel may vary with $\mu$.

The notion of convergence considered here is motivated by the spatial viewpoint of stochastic geometry, where random structures are studied through their restrictions to bounded regions of Euclidean space; see, for instance, \cite{penrose,last-penrose}. In particular, Lichev, Lodewijks, Mitsche and Schapira \cite[Definition~1.7]{lichev2023} introduced a closely related notion of spatial local convergence, based on convergence in distribution of graph restrictions to bounded Euclidean domains. Our approach follows the same principle of bounded spatial observations but requires convergence in total variation of the embedded graph distributions, retaining both vertex positions and their connections without imposing finite-range dependence. This differs from the classical Benjamini--Schramm local weak convergence, which compares rooted neighbourhoods in graph distance \cite{benjamini-schramm}. The resulting framework is particularly suited to studying the different limiting regimes of the GRIG.

Let $(G^{(n)})_{n\ge1}$ and $G^{(\infty)}$ be random spatial graphs whose underlying vertex processes have the same law, that of a homogeneous Poisson point process $V$ on $\R^d$ with intensity $\lambda>0$. For a bounded Borel set $W\subset\R^d$, let $G|_W$ denote the embedded graph induced by the vertices in $V\cap W$, retaining both their positions and their connections. Only the vertices are restricted to $W$; edges between them may still be supported by groups outside $W$.

\begin{definition}[Spatially local convergence in total variation]\label{def:spatial-local-TV}
We say that $G^{(n)}$ converges spatially locally in total variation to $G^{(\infty)}$ if, for every bounded Borel set $W\subset\R^d$,
\begin{equation}\label{eq:spatial-local-TV}
d_{\mathrm{TV}}\left(
  \mathcal L(G^{(n)}|_W),
  \mathcal L(G^{(\infty)}|_W)
\right)\longrightarrow0 
\qquad\text{as }n\to\infty,
\end{equation}
where $d_{\mathrm{TV}}(P,\tilde P):=\sup_A|P(A)-\tilde P(A)|$. The spatial window is kept fixed before taking the parameter limit.
\end{definition}

For a parametrized family, convergence along a prescribed parameter limit means that the definition holds along every sequence approaching that limit. This convention covers both $\mu\to\infty$ and $\mu\downarrow0$ below.

Since all the models considered below have the same Poisson vertex process, in the couplings used in the proofs we realize them on a common vertex configuration $V$.

\subsection{A total variation bound}\label{sec:tv-bound}

Conditionally on the vertex positions, dependencies between edges in the GRIG arise from groups joining three or more vertices. The following proposition bounds the total variation distance between the GRIG and a spatial graph with independent edges, including deterministic connection rules as a special case. For a measurable symmetric function $\varphi\colon\mathbb R^d\to[0,1]$, the random connection model $G^\varphi_{\mathrm{RCM}}$ has a homogeneous Poisson vertex process of intensity $\lambda$, and, conditionally on the positions, connects the vertices at $x$ and $y$ independently with probability $\varphi(x-y)$ \cite{meester-roy,last-nestmann-schulte}. When $\varphi$ takes values in $\{0,1\}$, it is a deterministic function of the vertex positions. The GRIG has connection function $p=1-e^{-\mu f}$ by \eqref{eq:edge-probability}, so that $G^p_{\mathrm{RCM}}$ has the same edge probabilities as $G$, but independent edges.

\begin{proposition}[Total variation bound]\label{prop:tv-bound}
Let $G$ be a GRIG with intensities $\lambda,\mu$ and membership kernel $g$ satisfying \eqref{eq:g-assumption}, and, for a bounded Borel set $W\subset\mathbb R^d$, let $\gamma_W(z):=\int_Wg(x-z)\dd x$. Then, for every measurable symmetric function $\varphi\colon\mathbb R^d\to[0,1]$,
\begin{equation}\label{eq:tv-bound}
    d_{\mathrm{TV}}\bigl(\mathcal L(G|_W),\mathcal L(G^\varphi_{\mathrm{RCM}}|_W)\bigr)\le \frac{2\lambda^3\mu}{3}\int_{\mathbb R^d}\gamma_W(z)^3\dd z +\frac{\lambda^2}{2}\iint_{W^2}\bigl|p(x-y)-\varphi(x-y)\bigr|\dd x\dd y . \end{equation}
If $\varphi$ takes values in $\{0,1\}$, the first term on the right-hand side can be omitted. \end{proposition}

\begin{proof}
We first couple $G$ with the random connection model $G^p_{\mathrm{RCM}}$, using the same Poisson vertex process $V$. Condition on $V\cap W=\{v_1,\ldots,v_k\}$ and write $\mathbf v=(v_1,\ldots,v_k)$. Throughout the following conditional estimates, probabilities are taken given this fixed vertex configuration. By Lemma~\ref{lem:types}, the counts $N_S$ of groups of exact type $S\subseteq[k]$, $|S|\ge2$, are independent Poisson variables with means $\Lambda_S(\mathbf v)$. In particular,
\[
\mu f(v_i-v_j)=\sum_{S\supseteq\{i,j\}}\Lambda_S(\mathbf v).
\]

To obtain independent edges with the same marginal probabilities, we retain the counts $N_{\{i,j\}}$ and replace the contributions of groups joining three or more vertices by independent Poisson variables. More precisely, for each pair $i<j$, introduce a Poisson variable $Z_{ij}$, independent of all the other variables, with conditional mean
\[
\mu f(v_i-v_j)-\Lambda_{\{i,j\}}(\mathbf v)
=\sum_{\substack{S\supseteq\{i,j\}\\|S|\ge3}}\Lambda_S(\mathbf v).
\]
The resulting counts $\widehat N_{ij}:=N_{\{i,j\}}+Z_{ij}$ are independent Poisson variables with means $\mu f(v_i-v_j)$. Connecting $v_i$ and $v_j$ whenever $\widehat N_{ij}\ge1$ therefore produces $G^p_{\mathrm{RCM}}$ on the same vertices.

If all counts $N_S$ with $|S|\ge3$ and all added variables $Z_{ij}$ vanish, both graphs have an edge between $v_i$ and $v_j$ precisely when $N_{\{i,j\}}\ge1$. Therefore, the graphs can differ only if at least one of these variables is positive. Since each type $S$ with $|S|\ge3$ contributes to the means of $\binom{|S|}{2}$ added variables, the union bound gives
\[
\Prob\bigl(G[\mathbf v]\neq G^p_{\mathrm{RCM}}[\mathbf v]\bigr)
\le\sum_{\substack{S\subseteq[k]\\|S|\ge3}}
\left(1+\binom{|S|}{2}\right)\Lambda_S(\mathbf v).
\]

We next couple $G^p_{\mathrm{RCM}}$ and $G^\varphi_{\mathrm{RCM}}$ using a common uniform variable for each pair of vertices. Conditionally on $\mathbf v$, their probability of disagreement is at most
\[
\sum_{i<j}\bigl|p(v_i-v_j)-\varphi(v_i-v_j)\bigr|.
\]

We now average the conditional error bounds over the vertex process. For a hypothetical group at position $z$, let $K_z$ denote the number of vertices in $V\cap W$ linked to it. By independent thinning of $V$, we have $K_z\sim\operatorname{Poisson}(\lambda\gamma_W(z))$. The Mecke formula gives
\[
\begin{aligned}
\E\Biggl[\sum_{\substack{S\subseteq[k]\\|S|\ge3}}\left(1+\binom{|S|}{2}\right)\Lambda_S(\mathbf V_W)\Biggr]
&=\mu\int_{\R^d}\E\left[\left(1+\binom{K_z}{2}\right)\1_{\{K_z\ge3\}}\right]\dd z\\
&\le4\mu\int_{\R^d}\E\left[\binom{K_z}{3}\right]\dd z =\frac{2\lambda^3\mu}{3}\int_{\R^d}\gamma_W(z)^3\dd z,
\end{aligned}
\]
where $\mathbf V_W$ is a measurable ordering of $V\cap W$, and we used $1+\binom{k}{2}\le4\binom{k}{3}$ for $k\ge3$, together with the third factorial moment of a Poisson variable.

By the multivariate Mecke formula, the expected disagreement bound between $G^p_{\mathrm{RCM}}$ and $G^\varphi_{\mathrm{RCM}}$ becomes
\[
\frac{\lambda^2}{2}\iint_{W^2}\bigl|p(x-y)-\varphi(x-y)\bigr|\dd x\dd y.
\]
Combining the two couplings and applying the coupling inequality proves \eqref{eq:tv-bound}.

Finally, suppose that $\varphi$ takes values in $\{0,1\}$. Conditionally on $\mathbf v$, the graph $G^\varphi_{\mathrm{RCM}}[\mathbf v]$ is deterministic, so we can compare it directly with $G[\mathbf v]$. Since each edge $\{v_i,v_j\}$ is present in the GRIG with probability $1-e^{-\mu f(v_i-v_j)}$, the union bound gives
\[
\Prob\bigl(G[\mathbf v]\neq G^\varphi_{\mathrm{RCM}}[\mathbf v]\bigr)
\le\sum_{i<j}\bigl|p(v_i-v_j)-\varphi(v_i-v_j)\bigr|.
\]
Averaging over the vertex process yields the second term in \eqref{eq:tv-bound}, but not the first.
\end{proof}

The bound in Proposition~\ref{prop:tv-bound} can be interpreted through the two-step approximation
\[
G\longrightarrow G^p_{\mathrm{RCM}}\longrightarrow G^\varphi_{\mathrm{RCM}}.
\]
The first step replaces the dependencies induced by groups joining multiple vertices with independent edges, while preserving the edge probabilities $p$; the second step replaces $p$ by the connection function $\varphi$ of the target graph. Accordingly, the first term in \eqref{eq:tv-bound} controls the dependence error: it is four times the bound $\frac{\lambda^3\mu}{6}\int\gamma_W(z)^3\dd z$ on the expected number of groups linked to at least three vertices in $W$. The second term is the distance between the two connection functions, integrated over the window. When $\varphi$ takes values in $\{0,1\}$, the target graph is deterministic given the vertex positions, the intermediate step is unnecessary, and only the second term remains.

\subsection{Scaling of the membership kernel}\label{sec:scaling}

For a GRIG with group intensity $\mu$ and membership kernel $g$ satisfying \eqref{eq:g-assumption}, define
\begin{equation}\label{eq:kappa}
\kappa_s(\mu,g):=\mu\int_{\R^d}g(z)^s\dd z=\mu\lVert g\rVert_s^s,\qquad s=1,2,\ldots.
\end{equation}
For convenience, we write $\kappa_s=\kappa_s(\mu,g)$, suppressing its dependence on $\mu$ and $g$,
whenever no confusion arises.
These quantities are finite since $0\le g\le1$ and $g$ is integrable. By the Mecke formula, $\kappa_1$ is the mean number of groups joined by a prescribed vertex. More generally, $s$ prescribed vertices at $v_1,\ldots,v_s$ share on average $\mu\int_{\R^d}\prod_{i=1}^s g(z-v_i)\dd z$ groups. By H\"older's inequality, this mean is at most $\kappa_s$, and by continuity of translations in $L^s(\R^d)$ it converges to $\kappa_s$ as the vertex positions coalesce. Thus, $\kappa_s$ is the supremum of the mean number of groups shared by $s$ prescribed vertices, and this supremum is approached as their pairwise distances tend to zero.

In particular, writing $f=g*g$ as in Section~\ref{sec:model}, we have $\kappa_2=\mu f(0)$. 
Since $f$ is continuous and, by the Cauchy--Schwarz inequality, maximal at the origin, 
the connection probability \eqref{eq:edge-probability} satisfies $p(x)\le1-e^{-\kappa_2}$, 
with $p(x)\to1-e^{-\kappa_2}$ as $x\to0$.

From now on, $g$ denotes a fixed reference kernel satisfying \eqref{eq:g-assumption}, and we consider families of GRIGs whose effective membership kernels are of the form
\begin{equation}\label{eq:two-scales}
h_{\theta,r}(x):=\theta g(x/r),\qquad \theta\in(0,1],\quad r>0.
\end{equation}
We call $\theta$ the \textit{strength} and $r$ the \textit{range} of the memberships. Hereafter, $f=g*g$ always denotes the overlap function of the reference kernel. For this model,
\[
\kappa_s(\mu,h_{\theta,r})=\mu\lVert h_{\theta,r}\rVert_s^s
=\mu r^d\theta^s\lVert g\rVert_s^s
=\kappa_s(\mu r^d,\theta g),\qquad s\ge1,
\]
so the quantities $\kappa_s$ are unchanged by the following spatial rescaling.
\begin{lemma}[Scaling]\label{lem:scaling}
Consider the GRIG with intensities $\lambda,\mu$ and membership kernel $h_{\theta,r}$ as in \eqref{eq:two-scales}. Under the spatial rescaling $x\mapsto x/r$ of both vertex and group positions, with all memberships preserved, the resulting configuration is distributed as the GRIG with intensities $\lambda r^d,\mu r^d$ and kernel $\theta g$. The same holds for the configuration with an additional vertex at the origin. Consequently, the projected graphs, rooted at the origin in the second case, have the same law up to graph isomorphism.
\end{lemma}

\begin{proof}
By the mapping theorem \cite[Theorem~5.1]{last-penrose}, the images of $V$ and $U$ under $x\mapsto x/r$ are independent homogeneous Poisson processes with intensities $\lambda r^d$ and $\mu r^d$. Conditionally on the positions, memberships remain independent, and a vertex at $v$ joins a group at $u$ with probability $\theta g((v-u)/r)=\theta g(v/r-u/r)$, which is precisely the membership probability of the rescaled configuration under the kernel $\theta g$. The origin is fixed by the rescaling, so the same argument applies to the configuration with an additional vertex at the origin. Since adjacency is defined through shared groups, the projected graphs are preserved up to isomorphism.
\end{proof}

The lemma shows that the range can be absorbed into the intensities, whereas the strength cannot in general be absorbed by a spatial rescaling. For the fixed reference kernel $g$, the parameters $(\mu r^d,\theta)$ of the rescaled model and the corresponding pair $(\kappa_1,\kappa_2)$ determine each other. Here and below, we abbreviate
\[
\kappa_s\equiv\kappa_s(\mu,h_{\theta,r})=\kappa_s(\mu r^d,\theta g).
\]
In particular,
\begin{equation}\label{eq:kappa-inverse}
\kappa_1=\mu r^d\theta\lVert g\rVert_1,\qquad
\kappa_2=\mu r^d\theta^2\lVert g\rVert_2^2
\quad\Longleftrightarrow\quad
\theta=\frac{\kappa_2\lVert g\rVert_1}{\kappa_1\lVert g\rVert_2^2},\qquad
\mu r^d=\frac{\kappa_1^2\lVert g\rVert_2^2}{\kappa_2\lVert g\rVert_1^2}.
\end{equation}
The constraint $\theta\le1$ is equivalent to $\kappa_2\le\kappa_1\lVert g\rVert_2^2/\lVert g\rVert_1$, and the higher-order quantities $\kappa_s=\mu r^d\theta^s\lVert g\rVert_s^s$ are also determined.

\paragraph{Interpretation.} The range determines the spatial scale of the memberships. By Lemma~\ref{lem:scaling}, its effect on observables invariant under spatial rescaling is captured by the rescaled intensities $\lambda r^d$ and $\mu r^d$. The strength controls how readily a group connects close vertices: the maximal membership probability is $\theta g(0)$, equal to $\theta$ when $g(0)=1$. Moreover, $\kappa_1$ is the mean number of groups joined by a vertex, while $\kappa_2/\kappa_1=\theta\lVert g\rVert_2^2/\lVert g\rVert_1$ is the limiting probability that a second close vertex joins a group sampled among those joined by the first vertex. Similarly, $\kappa_3/\kappa_2=\theta\lVert g\rVert_3^3/\lVert g\rVert_2^2$ describes the limiting probability that a third close vertex joins a typical group shared by the first two. These ratios therefore quantify the extent to which groups are shared among close vertices.

\subsection{Convergence to limiting random graph models}

We now identify the random graph models that arise as spatially local limits of the GRIG. Throughout this subsection, we use the parametrization \eqref{eq:two-scales}, with $\theta=\theta_\mu$ and $r=r_\mu$ depending on $\mu$.

\subsubsection{Limiting models}

All limiting graphs below have a homogeneous Poisson vertex process of intensity $\lambda$. While the RCM and RGG retain spatial dependence, the RIG, ER, and complete graph models have edge distributions that, conditional on the vertex process, do not depend on vertex positions.

\paragraph{Random connection models.}
For a measurable symmetric function $\varphi:\R^d\to[0,1]$, the random connection model $G_{\mathrm{RCM}}^\varphi$ connects vertices at $x$ and $y$ independently with probability $\varphi(x-y)$ \cite{meester-roy,last-nestmann-schulte}. These are the graphs $G^\varphi_{\mathrm{RCM}}$ of Section~\ref{sec:tv-bound}.

\paragraph{Random geometric graphs.} Suppose that $g$ has bounded support and let $R:=\sup\{\lVert x\rVert:g(x)>0\}$. Since $g$ is radial and non-increasing, $g>0$ on $B_R(0)$ and vanishes outside $\overline{B_R(0)}$. Consequently, $f(x)>0$ if and only if $\lVert x\rVert<2R$. The random geometric graph $G_{\mathrm{RGG}}$ connects two vertices whenever their distance is less than $2R$, and it thus corresponds to the random connection model with $\varphi=\1_{B_{2R}(0)}$.

\paragraph{Erd\H{o}s--R\'enyi and complete graphs.} For $p\in[0,1]$, the graph $G_{\mathrm{ER}}^p$ connects each pair of vertices independently with probability $p$, irrespective of their positions. Conditionally on $|V\cap W|=n$, its restriction to a bounded window $W$ is the Erd\H{o}s--R\'enyi graph $G(n,p)$. The complete graph $G_{\mathrm K}$ connects all pairs of vertices.

\paragraph{Random intersection graphs.} Let $U_\alpha$ be a homogeneous Poisson point process on $\R^d$ with intensity $\alpha>0$, independent of the vertices. Conditionally on $U_\alpha$, each vertex independently joins each group $u\in U_\alpha$ with probability $g(u)$, regardless of its spatial position, and two vertices are adjacent whenever they share a group. We denote this graph by $G_{\mathrm{RIG}}^{(\alpha,g)}$.

Since vertex positions play no role, the law of this graph depends on $(\alpha,g)$ only through the image measure $\nu$ of $\alpha\,\mathrm du$ under $u\mapsto g(u)$, restricted to $(0,1]$. Indeed, by the mapping theorem, the positive membership probabilities $\{g(u):u\in U_\alpha,\ g(u)>0\}$ form a Poisson process on $(0,1]$ with intensity measure $\nu$, satisfying $\nu([\varepsilon,1])\le\alpha\lVert g\rVert_1/\varepsilon$. In particular, $G_{\mathrm{RIG}}^{(\alpha,g(\cdot/r))}$ and $G_{\mathrm{RIG}}^{(\alpha r^d,g)}$ have the same law. Although the total number of groups may be infinite, each vertex joins a Poisson number of groups with mean $\alpha\lVert g\rVert_1$, so the restriction to a bounded window involves finitely many groups almost surely. By the probability generating functional, two vertices are adjacent with probability $1-e^{-\alpha\lVert g\rVert_2^2}$.

\subsubsection{Convergence}
For the membership kernel $h_{\theta,r}=\theta g(\cdot/r)$, a change of variables gives the mean number of groups shared by two prescribed vertices at displacement $x$:
\begin{equation}\label{eq:pair-intensity}
\mu\theta^2 r^d f(x/r)
=\frac{\kappa_2}{\lVert g\rVert_2^2}f(x/r),
\qquad
\kappa_2=\mu\theta^2 r^d\lVert g\rVert_2^2.
\end{equation}
Here $f$ is continuous and $f\le f(0)=\lVert g\rVert_2^2$, as observed in Section~\ref{sec:scaling}.
Let us also recall that the dependence of $\theta$ and $r$ on $\mu$ is left implicit.

\begin{theorem}[Spatially local limits]\label{thm:local-limits}
Under \eqref{eq:g-assumption}, let $G^{(\mu)}$ be the GRIG with fixed vertex intensity $\lambda$, group intensity $\mu$, and membership kernel
\(
h_{\theta,r},\,\theta\in(0,1],\, r>0.
\)
For a bounded Borel set $W\subset\R^d$ and a random spatial graph $G^{(\infty)}$, write
\[
d_{\mu,W}(G^{(\infty)}):=d_{\mathrm{TV}}\bigl(\mathcal L(G^{(\mu)}|_W),\mathcal L(G^{(\infty)}|_W)\bigr),
\qquad
\gamma_W(z):=\int_W g(x-z)\dd x.
\]
In each of the following cases, $G^{(\mu)}$ converges spatially locally in total variation to the indicated graph, with the stated bound. We let $\mu\to\infty$ in cases (i), (ii), (iv), and (v), and $\mu\downarrow0$ in case (iii).
\begin{enumerate}
\item[(i)] \emph{(RCM)} If $r=1$ and $\mu\theta^2\to c\in(0,\infty)$, the limit is $G_{\mathrm{RCM}}^\varphi$, with $\varphi=1-e^{-cf}$, and
\[
d_{\mu,W}(G_{\mathrm{RCM}}^\varphi)
\le\frac{2\lambda^3}{3}\mu\theta^3\int_{\R^d}\gamma_W(z)^3\dd z
+\frac{\lambda^2|W|^2\lVert g\rVert_2^2}{2}|\mu\theta^2-c|
\underset{\mu\to\infty}{\longrightarrow}0,
\]
where $\int_{\R^d}\gamma_W(z)^3\dd z\le|W|\lVert g\rVert_1^3$.

\item[(ii)] \emph{(RGG)} If $r=1$, $\mu\theta^2\to\infty$, and $g$ has bounded support, the limit is $G_{\mathrm{RGG}}$, with connection radius $2R$, and
\[
d_{\mu,W}(G_{\mathrm{RGG}})
\le\frac{\lambda^2}{2}\int_{W^2}
\1_{\{\lVert x-y\rVert<2R\}}
e^{-\mu\theta^2 f(x-y)}\dd x\dd y
\underset{\mu\to\infty}{\longrightarrow}0.
\]

\item[(iii)] \emph{(RIG)} If $\theta=1$ and $r=(\alpha/\mu)^{1/d}$ for a fixed $\alpha>0$, the limit is $G_{\mathrm{RIG}}^{(\alpha,g)}$, and, writing $\omega_g(h):=\lVert g(\cdot-h)-g\rVert_1$,
\[
d_{\mu,W}(G_{\mathrm{RIG}}^{(\alpha,g)})
\le\alpha\lambda\int_W\omega_g(x/r)\dd x
\underset{\mu\downarrow0}{\longrightarrow}0.
\]

\item[(iv)] \emph{(ER)} If $r\to\infty$, $\theta\to0$, and $\mu\theta^2r^d\to c\in(0,\infty)$, the limit is $G_{\mathrm{ER}}^p$, with $p=1-e^{-c\lVert g\rVert_2^2}$, and
\[
\begin{aligned}
d_{\mu,W}(G_{\mathrm{ER}}^p)
 &\le\frac{2\lambda^3|W|^3\lVert g\rVert_1}{3\lVert g\rVert_2^2}\,\theta\kappa_2(\mu, h_{\theta,r})\\
&\quad+\frac{\lambda^2}{2}\int_{W^2}
\left|e^{-\kappa_2(\mu, h_{\theta,r}) f((x-y)/r)/\lVert g\rVert_2^2}
-e^{-c\lVert g\rVert_2^2}\right|\dd x\dd y
\underset{\mu\to\infty}{\longrightarrow}0.
\end{aligned}
\]

\item[(v)] \emph{(Complete graph)} If $r\to\infty$ and $\mu\theta^2r^d\to\infty$, the limit is $G_{\mathrm K}$, and
\[
d_{\mu,W}(G_{\mathrm K})
\le\frac{\lambda^2}{2}\int_{W^2}
e^{-\kappa_2(\mu, h_{\theta,r}) f((x-y)/r)/\lVert g\rVert_2^2}
\dd x\dd y
\underset{\mu\to\infty}{\longrightarrow}0.
\]
\end{enumerate}
\end{theorem}

The five rows of Table~\ref{tab:regimes} are particular choices satisfying respectively cases (i)--(v) of the theorem.

\begin{proof}
In cases \emph{(i), (ii), (iv),} and \emph{(v)}, we apply Proposition~\ref{prop:tv-bound} with membership kernel $h_{\theta,r}$. Its pairwise intensity is given by \eqref{eq:pair-intensity}, while the corresponding window function is $z\mapsto\theta\int_W g((x-z)/r)\dd x$.

\emph{(i)} Take $\varphi(x)=1-e^{-cf(x)}$. Since $r=1$, the window function is $\theta\gamma_W$. Moreover,
\[
\int_{\R^d}\gamma_W(z)^3\dd z
\le\lVert\gamma_W\rVert_\infty^2\lVert\gamma_W\rVert_1
\le|W|\lVert g\rVert_1^3.
\]
Using $|e^{-u}-e^{-v}|\le|u-v|$ for $u,v\ge0$ and $f\le\lVert g\rVert_2^2$, the second term of \eqref{eq:tv-bound} is bounded by $\frac{\lambda^2|W|^2\lVert g\rVert_2^2}{2}|\mu\theta^2-c|$. The first term vanishes because $\mu\theta^2\to c$ implies $\theta\to0$ and hence $\mu\theta^3=\theta(\mu\theta^2)\to0$, as $\mu\to\infty$.

\emph{(ii)} Take $\varphi(x)=\1_{\{\lVert x\rVert<2R\}}$, so that $G_{\mathrm{RCM}}^\varphi=G_{\mathrm{RGG}}$. The deterministic case of Proposition~\ref{prop:tv-bound} applies. If $\lVert x-y\rVert\ge2R$, then $f(x-y)=0$ and the two edge probabilities agree. Otherwise, the disagreement probability is $e^{-\mu\theta^2f(x-y)}$, which tends to zero because $f>0$ on $B_{2R}(0)$ and $\mu\theta^2\to\infty$. Dominated convergence on $W^2$ concludes.

\emph{(iii)} Rescale the group positions in the GRIG by writing $u=z/r$. Since $\mu r^d=\alpha$, the rescaled groups form a homogeneous Poisson point process of intensity $\alpha$, which we use as the common group process $U_\alpha$ for both models. Conditionally on the vertex positions and on $U_\alpha$, a vertex at $v$ joins a group at $u$ with probability $g(u-v/r)$ in the GRIG and $g(u)$ in the limiting RIG. For every vertex--group pair, we use the same independent uniform variable to generate the membership indicators, which then differ with probability $|g(u-v/r)-g(u)|$. By the Mecke formula, the expected number of differing memberships involving vertices in $W$ is
\[
\alpha\lambda\int_W\int_{\R^d}|g(u-v/r)-g(u)|\dd u\dd v
=\alpha\lambda\int_W\omega_g(v/r)\dd v.
\]
If no membership differs, the induced graphs on $V\cap W$ coincide, including their vertex positions. The coupling inequality and Markov's inequality give the desired bound. Finally, by assumption $r\to\infty$ as $\mu\downarrow0$, translations are continuous in $L^1(\R^d)$, and $\omega_g\le2\lVert g\rVert_1$, so dominated convergence concludes.

\emph{(iv)} Take $\varphi\equiv1-e^{-c\lVert g\rVert_2^2}=p$, so that $G_{\mathrm{RCM}}^\varphi=G_{\mathrm{ER}}^p$. The window function satisfies
\[
\theta\int_W g((x-z)/r)\dd x\le\theta|W|,
\]
and its integral over $z\in\R^d$ equals $\theta r^d|W|\lVert g\rVert_1$. Consequently, the first term of \eqref{eq:tv-bound} is at most
\[
\frac{2\lambda^3}{3}\mu\theta^3r^d|W|^3\lVert g\rVert_1
=\frac{2\lambda^3|W|^3\lVert g\rVert_1}{3\lVert g\rVert_2^2}\,\theta\kappa_2(\mu, h_{\theta,r})
\underset{\mu\to\infty}{\longrightarrow}0,
\]
because $\mu\theta^2r^d\to c$ and $\theta\to0$ as $\mu\to\infty$. For the second term, continuity of $f$ gives $f((x-y)/r)\to f(0)=\lVert g\rVert_2^2$, while $\kappa_2(\mu, h_{\theta,r})\to c\lVert g\rVert_2^2$. The integrand therefore tends to zero and is bounded by one, so dominated convergence applies.

\emph{(v)} Take $\varphi\equiv1$, so that $G_{\mathrm{RCM}}^\varphi=G_{\mathrm K}$. The deterministic case of Proposition~\ref{prop:tv-bound} gives the stated bound without a dependence term. Since $r\to\infty$, we have $f((x-y)/r)\to f(0)=\lVert g\rVert_2^2>0$, while $\kappa_2(\mu, h_{\theta,r})=\mu\theta^2r^d\lVert g\rVert_2^2\to\infty$. The integrand tends to zero and is bounded by one, so dominated convergence concludes.
\end{proof}

Table~\ref{tab:regimes} lists five representative scalings and the corresponding values of $\kappa_1$ and $\kappa_2$.

\begin{table}[ht]
\centering
\small
\renewcommand{\arraystretch}{1.35}
\begin{tabular*}{0.8\textwidth}{@{\extracolsep{\fill}}c|cc|cc@{}}
\toprule
Regime & $\theta$ & $r$ & $\kappa_1$ & $\kappa_2$\\
\midrule
RCM ($\mu \to \infty$)& $\mu^{-1/2}$ & $1$ & $\mu^{1/2}\lVert g\rVert_1$ & $\lVert g\rVert_2^2$ \\
RGG ($\mu \to \infty$)& $\mu^{-a}, a\in[0,\frac12)$ & $1$ & $\mu^{1-a}\lVert g\rVert_1$ & $\mu^{1-2a}\lVert g\rVert_2^2$ \\
RIG ($\mu \downarrow 0$)& $1$ & $(\alpha/\mu)^{1/d}$ & $\alpha\lVert g\rVert_1$ & $\alpha\lVert g\rVert_2^2$ \\
ER ($\mu \to \infty$) & $\mu^{-1}$ & $\mu^{1/d}$ & $\mu\lVert g\rVert_1$ & $\lVert g\rVert_2^2$ \\
Complete ($\mu \to \infty$) & $1$ & $\mu^{1/d}$ & $\mu^2\lVert g\rVert_1$ & $\mu^2\lVert g\rVert_2^2$ \\
\bottomrule
\end{tabular*}
\caption{Representative parameter choices for the five limiting models of Theorem~\ref{thm:local-limits}, with membership kernels of the form $h_{\theta,r}=\theta g(\cdot/r)$ and $\kappa_s=\mu\theta^s r^d\lVert g\rVert_s^s$. Here $\alpha>0$ is fixed.}
\label{tab:regimes}
\end{table}

\begin{remark}[Groups joining three vertices]\label{rem:triple-groups}
For $r=1$, the Mecke formula and the third factorial moment of the number of vertices linked to a prescribed group give
\[
\begin{aligned}
&\Prob\left(\exists\text{ a group joined to at least three vertices in }V\cap W\right)\\
&\qquad\le\frac{\lambda^3}{6}\mu\theta^3
\int_{\R^d}\gamma_W(z)^3\dd z
=\frac{\lambda^3\kappa_3}{6\lVert g\rVert_3^3}
\int_{\R^d}\gamma_W(z)^3\dd z.
\end{aligned}
\]
Such groups therefore disappear whenever $\kappa_3\to0$, as in case (i). Under the RGG scaling $\theta=\mu^{-a}$, we have $\kappa_3=\mu^{1-3a}\lVert g\rVert_3^3$, so the bound shows that they disappear for $a\in(1/3,1/2)$. For $a\le1/3$, the bound does not tend to zero and gives no conclusion. This distinction is not needed for the RGG convergence in case (ii), which holds throughout $a\in[0,1/2)$. The threshold $a=1/3$ reappears in the clique--cycle decomposition of clustering; see Theorem~\ref{thm:rgg-clique-cycle}.
\end{remark}

\begin{remark}[Poisson mixture of binomial RIGs]
For the hard-threshold kernel $g=p\1_{B_R(0)}$, with $p\in(0,1]$, the RIG limit has a particularly simple interpretation. Only groups in $B_R(0)$ have a positive membership probability, and their number $m$ is Poisson distributed with mean $\alpha|B_R(0)|$. Conditionally on $|V\cap W|=n$ and on $m$, the restriction of the limiting graph to $W$ has the law of the classical binomial random intersection graph $G(n,m,p)$. Thus, conditionally on $|V\cap W|=n$, the limit is a Poisson mixture of binomial RIGs. For a general kernel, the limiting construction allows heterogeneous membership probabilities, described by the measure $\nu$.
\end{remark}

\begin{remark}[Local finiteness]\label{rem:local-finiteness}
In the GRIG, every vertex has finite degree almost surely under \eqref{eq:g-assumption}. In contrast, in $G_{\mathrm{RIG}}^{(\alpha,g)}$, every vertex belonging to at least one group has infinitely many neighbours almost surely, since each such group is joined independently by infinitely many vertices of the underlying Poisson process. The same holds for $G_{\mathrm{ER}}^p$ with $p>0$ and for $G_{\mathrm K}$. Accordingly, in cases (iii)--(v), where $r\to\infty$, the mean degree of the GRIG satisfies
\[
\E[D(0)]
=\lambda\int_{\R^d}\left(1-e^{-\mu\theta^2r^df(x/r)}\right)\dd x
=\lambda r^d\int_{\R^d}
\left(1-e^{-\kappa_2f(x)/\lVert g\rVert_2^2}\right)\dd x
\longrightarrow\infty.
\]
Nevertheless, the restrictions to bounded windows contain finitely many vertices almost surely, and Theorem~\ref{thm:local-limits} is well defined. This illustrates why spatially local convergence is appropriate here, whereas convergence based on finite graph-distance neighbourhoods would require additional care.
\end{remark}

\section{Clustering in asymptotic regimes}\label{sec:clustering-limits}

In Section~\ref{sec:clustering}, we obtained explicit representations for the global and average local clustering coefficients of the GRIG, while Section~\ref{sec:graph-limits} established its convergence to several limiting random graph models under suitable parameter scalings. We now investigate how the clustering coefficients behave in these asymptotic regimes, and whether their limits agree with the clustering of the corresponding limiting models.

A central difficulty is that spatially local convergence does not necessarily imply clustering convergence. Indeed, the latter depends on the structure of the neighbourhood of a typical vertex, which may involve connections at arbitrarily large spatial distances. As a result, contributions from increasingly distant vertices may remain relevant even when the graph converges on every fixed bounded window. Understanding when clustering is preserved by the limiting procedure, and when it is not, is one of the main objectives of this section.

We show that, for a fixed reference kernel, global clustering depends on the scaling parameters only through $\kappa_1$ and $\kappa_2$. We then analyse its boundary behaviour and apply the results to the five benchmark regimes of Table \ref{tab:regimes}. Finally, we explain why the local limit and the clustering limit agree in some regimes and not in others, and we study the clique--cycle decomposition under saturation.

\subsection{Clustering as a function of two parameters}\label{sec:two-parameters}

As an application of Lemma~\ref{lem:scaling}, we first show that the range of the membership kernel affects the conditional quantities of Section~\ref{sec:clustering} only through a spatial rescaling of the intensities.

\begin{corollary}[Scaling of clustering]\label{cor:clustering-scaling}
Let $\theta\in(0,1]$ and $r>0$. Let $(M,Q,Q_{\mathrm{cl}},Q_{\mathrm{cyc}})$ be the conditional quantities for intensities $\lambda,\mu$ and kernel $\theta g(\cdot/r)$. Denote by $(M',Q',Q'_{\mathrm{cl}},Q'_{\mathrm{cyc}})$ their counterparts for intensities $\lambda r^d,\mu r^d$ and kernel $\theta g$. Then
\[
(M,Q,Q_{\mathrm{cl}},Q_{\mathrm{cyc}})\overset{d}{=}(r^dM',Q',Q'_{\mathrm{cl}},Q'_{\mathrm{cyc}}).
\]
Consequently, the two models have the same global clustering coefficient and clique and cycle contributions, as well as the same mean degree. Their average local clustering coefficients also coincide, and the same holds for the degree-conditioned clustering spectrum $\bar c_\lambda(k)$, with the vertex intensity of the rescaled model equal to $\lambda r^d$.
\end{corollary}

\begin{proof}
By Lemma~\ref{lem:scaling}, we may couple the two models so that the map $x\mapsto x/r$ sends the groups $U_0$ joined by the origin to the corresponding groups $U_0'$ of the rescaled model. On this coupling,
\[
q(x)=1-\prod_{u\in U_0}\bigl(1-\theta g((x-u)/r)\bigr)=q'(x/r),
\]
hence, by \eqref{eq:M}, $M=r^dM'$. Similarly, \eqref{eq:s_e_eta} gives $s(x,y)=s'(x/r,y/r)$ and
\[
\eta(x,y)=\mu\int_{\R^d}\bigl(1-\theta g(z/r)\bigr)\theta g((z-x)/r)\theta g((z-y)/r)\dd z=\eta'(x/r,y/r),
\]
where we used the substitution $z=rz'$ and the fact that $\eta'$ is computed with group intensity $\mu r^d$. Substituting $x=rx'$ and $y=ry'$ in the integrals of Lemma~\ref{lem:Q-integral} yields the stated identity in law. The remaining assertions follow from \eqref{eq:global-group-representation}, \eqref{eq:local-global} and \eqref{eq:degree-clustering}.
\end{proof}

The scaling corollary shows that the range can be absorbed into the intensities without changing the global clustering coefficient, whereas average local clustering also depends on the rescaled vertex intensity $\lambda r^d$. We now show that, for a fixed reference kernel $g$, the remaining dependence of global clustering is captured by the two parameters $(\kappa_1,\kappa_2)$ of Section~4.2.

For comparison, let $\varphi:\R^d\to[0,1]$ be measurable and symmetric, with $0<\int_{\R^d}\varphi<\infty$. The random connection model $G_{\mathrm{RCM}}^\varphi$ has independent edges conditional on the vertex positions, and its global and average local clustering coefficients coincide and are given by
\begin{equation}\label{eq:rcm-clustering}
C_{\mathrm{RCM}}(\varphi):=\frac{\displaystyle\int_{\R^{2d}}\varphi(x)\varphi(y)\varphi(x-y)\dd x\dd y}{\displaystyle\left(\int_{\R^d}\varphi(x)\dd x\right)^2}.
\end{equation}
If $g$ has bounded support, let $R$ be as in Section~\ref{sec:graph-limits} and $B:=B_{2R}(0)$. The clustering coefficient of $G_{\mathrm{RGG}}$ is
\begin{equation}\label{eq:rgg-clustering}
C_{\mathrm{RGG}}:=C_{\mathrm{RCM}}(\1_B)=\frac{1}{|B|^2}\int_{B\times B}\1_{\{\lVert x-y\rVert<2R\}}\dd x\dd y,
\end{equation}
which, after rescaling by $2R$, depends only on the dimension $d$.

For the fixed reference kernel $g$, recall $f=g*g$ and define the three-point overlap function
\begin{equation}\label{eq:three-point-overlap}
A(x,y):=\int_{\R^d}g(z)g(z-x)g(z-y)\dd z.
\end{equation}
Since $0\le g\le1$, we have $A(x,y)\le\min\{f(x),f(y),f(x-y)\}$. By \eqref{eq:kappa-inverse}, the admissible parameter domain is
\[
\mathcal D_g:=\left\{(u,v)\in(0,\infty)^2:v\le u\frac{\lVert g\rVert_2^2}{\lVert g\rVert_1}\right\}.
\]
Each $(u,v)\in\mathcal D_g$ uniquely determines the strength $\theta$ and the group intensity $\nu$ of a unit-range GRIG with membership kernel $\theta g$, through $\kappa_1=u$ and $\kappa_2=v$. In particular, its strength is $\theta=v\lVert g\rVert_1/(u\lVert g\rVert_2^2)\le1$. Define
\[
a_v(z):=\frac{v}{\lVert g\rVert_2^2}f(z),\qquad \tau_{u,v}(x,y):=\frac{v^2\lVert g\rVert_1}{u\lVert g\rVert_2^4}A(x,y).
\]
These are the mean numbers of groups shared, respectively, by two vertices at displacement $z$ and by three vertices at positions $0,x,y$ in that unit-range model. In particular,
\[
0\le\tau_{u,v}(x,y)\le\min\{a_v(x),a_v(y),a_v(x-y)\}.
\]

Recall the configuration integrals $I_\wedge$, $I_{\mathrm{cl}}$, $I_{\mathrm{cyc}}$ and $I_\triangle$ from (3.7)--(3.8). For $(u,v)\in\mathcal D_g$, we denote by $I_i(u,v)$, $i\in\{\wedge,\mathrm{cl},\mathrm{cyc},\triangle\}$, the corresponding integrals evaluated in the unit-range GRIG determined by $(u,v)$. In terms of $a_v$ and $\tau_{u,v}$, their explicit expressions are
\begin{equation}\label{eq:two-param-wedge}
I_\wedge(u,v)=\int_{\R^{2d}}\Big[(1-e^{-a_v(x)})(1-e^{-a_v(y)})+e^{-a_v(x)-a_v(y)}(e^{\tau_{u,v}(x,y)}-1)\Big]\dd x\dd y,
\end{equation}
and
\begin{equation}\label{eq:two-param-triangle}
I_\triangle(u,v)=\int_{\R^{2d}}\Bigg[1-e^{-\tau_{u,v}(x,y)}+e^{-\tau_{u,v}(x,y)}\prod_{z\in\{x,y,x-y\}}\left(1-e^{-(a_v(z)-\tau_{u,v}(x,y))}\right)\Bigg]\dd x\dd y.
\end{equation}
The clique and cycle integrals are, correspondingly,
\[
I_{\mathrm{cl}}(u,v):=\int_{\R^{2d}}\bigl(1-e^{-\tau_{u,v}(x,y)}\bigr)\dd x\dd y,\qquad I_{\mathrm{cyc}}(u,v):=I_\triangle(u,v)-I_{\mathrm{cl}}(u,v).
\]

\begin{proposition}[Two-parameter representation]\label{prop:two-parameter}
Under \eqref{eq:g-assumption}, for any GRIG with intensities $\lambda,\mu>0$ and membership kernel $h_{\theta,r}$ as in \eqref{eq:two-scales}, the global clustering coefficient and its clique and cycle contributions satisfy
\begin{equation}\label{eq:two-param-representation}
C_{\mathrm{glob}}=\frac{I_\triangle(\kappa_1,\kappa_2)}{I_\wedge(\kappa_1,\kappa_2)},\qquad C_{\mathrm{glob}}^{(i)}=\frac{I_i(\kappa_1,\kappa_2)}{I_\wedge(\kappa_1,\kappa_2)},\quad i\in\{\mathrm{cl},\mathrm{cyc}\}.
\end{equation}
In particular, for fixed $g$, these coefficients depend only on $(\kappa_1,\kappa_2)$, independently of the vertex intensity $\lambda$ and the membership range $r$.
\end{proposition}

\begin{proof}
By Corollary~\ref{cor:clustering-scaling}, it suffices to consider the unit-range GRIG with kernel $\theta g$ and group intensity $\nu:=\mu r^d$. By \eqref{eq:kappa-inverse}, $(\nu,\theta)$ is uniquely determined by $(\kappa_1,\kappa_2)$. For three vertices at $x,0,y$, the mean number of groups shared by two vertices at displacement $z$ is
\[
\nu\theta^2f(z)=a_{\kappa_2}(z),
\]
whereas the mean number of groups joining all three vertices is
\[
\nu\theta^3A(x,y)=\tau_{\kappa_1,\kappa_2}(x,y).
\]
Proposition~3.1 therefore gives the integrands in \eqref{eq:two-param-wedge} and \eqref{eq:two-param-triangle}, as well as the corresponding clique and cycle terms. For the original kernel $h_{\theta,r}$, the change of variables $x=rx'$ and $y=ry'$ shows that each configuration integral of Section~\ref{sec:clustering} acquires the same factor $r^{2d}$ relative to its unit-range counterpart. This factor cancels in all the clustering ratios, and Corollary~3.3 gives \eqref{eq:two-param-representation}. Positivity of the denominator follows from Proposition~3.2.
\end{proof}

From now on we adopt again the notation $\theta=\theta_\mu$ and $r=r_\mu$. For a family with group intensity $\mu$, kernel $h_{\theta,r}=\theta g(\cdot/r)$ and vertex intensity $\lambda$, we write $C_{\mathrm{glob},\mu}$ for the global clustering coefficient and $C_{\mathrm{loc},\mu}(\lambda)$ for average local clustering evaluated at vertex intensity $\lambda$. The dependence of $\kappa_1,\kappa_2$ on $\mu$, $\theta$ and $r$ is left implicit.

The domain $\mathcal D_g$ has three asymptotic boundary regimes. Indeed, if $\kappa_2$ stays in a compact subset of $(0,\infty)$ and $\kappa_1$ stays bounded, then $(\kappa_1,\kappa_2)$ stays in a compact subset of $\mathcal D_g$, since $\kappa_1\ge\kappa_2\lVert g\rVert_1/\lVert g\rVert_2^2$. Thus any sequence leaving all compact subsets of $\mathcal D_g$ has a subsequence along which $\kappa_2\to0$, or $\kappa_2\to\infty$, or $\kappa_1\to\infty$ with $\kappa_2$ converging to a positive finite limit. The next theorem identifies the clustering limits in these regimes.

\begin{theorem}[Boundary behaviour of clustering]\label{thm:boundary}
Under \eqref{eq:g-assumption}, fix the reference kernel $g$ and consider a family of GRIGs with group intensities $\mu$, kernels $h_{\theta,r}=\theta g(\cdot/r)$, where $\theta\in(0,1]$ and $r>0$, and arbitrary vertex intensities $\lambda>0$. Let $\kappa_1,\kappa_2$ denote the corresponding parameters from \eqref{eq:kappa-inverse}. Then the following hold along any limiting regime of the family.
\begin{enumerate}
\item[(a)] \emph{(Weak sharing)} For every parameter choice with $\kappa_2\le1/4$,
\[
\left|C_{\mathrm{glob},\mu}-\frac{1}{1+\kappa_1}\right|\le6\kappa_2,\qquad C_{\mathrm{glob},\mu}^{(\mathrm{cyc})}\le2\kappa_2.
\]
\item[(b)] \emph{(Many groups)} If $\kappa_1\to\infty$ and $\kappa_2\to\kappa\in(0,\infty)$, then $C_{\mathrm{glob},\mu}\to C_{\mathrm{RCM}}(\varphi)$, where $\varphi=1-e^{-\kappa f/\lVert g\rVert_2^2}$.
\item[(c)] \emph{(Saturation)} If $g$ has bounded support and $\kappa_2\to\infty$, then $C_{\mathrm{glob},\mu}\to C_{\mathrm{RGG}}$.
\end{enumerate}
Moreover, in cases (b) and (c), $C_{\mathrm{loc},\mu}(\lambda)$ converges to $C_{\mathrm{RCM}}(\varphi)$ and $C_{\mathrm{RGG}}$, respectively, for every choice of the vertex intensities $\lambda$.
\end{theorem}

\begin{proof}
By Corollary~\ref{cor:clustering-scaling}, it suffices to work with the unit-range GRIG having group intensity $\nu:=\mu r^d$, vertex intensity $\lambda r^d$ and kernel $\theta g$. For brevity, set
\[
a_z:=\nu\theta^2f(z)=a_{\kappa_2}(z),\qquad \tau_{x,y}:=\nu\theta^3A(x,y)=\tau_{\kappa_1,\kappa_2}(x,y).
\]
We also have $\kappa_1=\nu\theta\lVert g\rVert_1$, $\kappa_2=\nu\theta^2\lVert g\rVert_2^2$ and $\kappa_3=\nu\theta^3\lVert g\rVert_3^3$. Throughout the proof, we abbreviate $I_\wedge:=I_\wedge(\kappa_1,\kappa_2)$, and similarly for $I_\triangle$, $I_{\mathrm{cl}}$ and $I_{\mathrm{cyc}}$. We use the identities
\begin{equation}\label{eq:f-A-integrals}
\int_{\R^d}f(x)\dd x=\lVert g\rVert_1^2,\qquad \int_{\R^d}A(x,y)\dd y=\lVert g\rVert_1f(x),\qquad \int_{\R^{2d}}A(x,y)\dd x\dd y=\lVert g\rVert_1^3,
\end{equation}
which follow from Tonelli's theorem, as well as the bounds $f\le\lVert g\rVert_2^2$, $A\le\lVert g\rVert_2^2$ and $\lVert f\rVert_2^2\le\lVert f\rVert_\infty\lVert f\rVert_1\le\lVert g\rVert_2^2\lVert g\rVert_1^2$.

\emph{(a) Weak sharing.} For $u,v\ge0$ and $0\le w\le\min\{u,v\}$, we have
\begin{equation}\label{eq:elementary}
\begin{aligned}
\bigl|(1-e^{-u})(1-e^{-v})-uv\bigr|&\le\tfrac12(u^2v+uv^2),\\
\bigl|e^{-u-v}(e^w-1)-w\bigr|&\le w(u+v),\\
\bigl|1-e^{-w}-w\bigr|&\le\tfrac12w^2.
\end{aligned}
\end{equation}
The first inequality follows from $\max\{0,u-u^2/2\}\le1-e^{-u}\le u$, the lower bound being trivial when $\max\{u,v\}>2$. For the second, use $we^{-u-v}\le e^{-u-v}(e^w-1)\le we^{w-u-v}\le w$ and $e^{-u-v}\ge1-u-v$; the third follows from Taylor's formula. Moreover, the product in \eqref{eq:two-param-triangle} lies in $[0,a_xa_ya_{x-y}]$. Set
\[
W_0:=(\nu\theta^2)^2\lVert g\rVert_1^4+\nu\theta^3\lVert g\rVert_1^3,\qquad T_0:=\nu\theta^3\lVert g\rVert_1^3.
\]
By \eqref{eq:f-A-integrals}, $\int\!\!\int(a_xa_y+\tau_{x,y})=W_0$ and $\int\!\!\int\tau_{x,y}=T_0$. Thus $I_\wedge=W_0+E_\wedge$ and $I_\triangle=T_0+E_\triangle$, where
\begin{align*}
|E_\wedge|&\le\int_{\R^{2d}}\Bigl[\tfrac12(a_x^2a_y+a_xa_y^2)+\tau_{x,y}(a_x+a_y)\Bigr]\dd x\dd y\\
&=(\nu\theta^2)^3\lVert f\rVert_2^2\lVert g\rVert_1^2+2\nu\theta^3\nu\theta^2\lVert g\rVert_1\lVert f\rVert_2^2\le2\kappa_2W_0,\\
|E_\triangle|&\le\int_{\R^{2d}}\Bigl[\tfrac12\tau_{x,y}^2+a_xa_ya_{x-y}\Bigr]\dd x\dd y\\
&\le\tfrac12\nu\theta^3\lVert g\rVert_2^2\,T_0+\nu\theta^2\lVert g\rVert_2^2(\nu\theta^2)^2\lVert g\rVert_1^4\le\kappa_2W_0.
\end{align*}
Here we used $\tau_{x,y}\le\nu\theta^3\lVert g\rVert_2^2\le\kappa_2$ and $f(x-y)\le\lVert g\rVert_2^2$. Moreover,
\[
\frac{T_0}{W_0}=\frac{1}{1+\nu\theta\lVert g\rVert_1}=\frac{1}{1+\kappa_1}=: \rho.
\]
For $\kappa_2\le1/4$, we therefore obtain
\[
\left|C_{\mathrm{glob},\mu}-\rho\right|=\left|\frac{W_0E_\triangle-T_0E_\wedge}{W_0(W_0+E_\wedge)}\right|\le\frac{\kappa_2+2\rho\kappa_2}{1-2\kappa_2}\le6\kappa_2.
\]
Finally, by \eqref{eq:two-param-triangle},
\[
I_{\mathrm{cyc}}\le\int_{\R^{2d}}a_xa_ya_{x-y}\dd x\dd y\le\kappa_2W_0,
\]
so $C_{\mathrm{glob},\mu}^{(\mathrm{cyc})}\le\kappa_2/(1-2\kappa_2)\le2\kappa_2$.

\emph{(b) Many groups.} By \eqref{eq:kappa-inverse}, $\theta=\kappa_2\lVert g\rVert_1/(\kappa_1\lVert g\rVert_2^2)\to0$ and $\nu\theta^2=\kappa_2/\lVert g\rVert_2^2\to\kappa/\lVert g\rVert_2^2$. Since $0\le\tau_{x,y}\le\theta\kappa_2\to0$ uniformly, the integrands of \eqref{eq:two-param-wedge} and \eqref{eq:two-param-triangle} converge pointwise to $\varphi(x)\varphi(y)$ and $\varphi(x)\varphi(y)\varphi(x-y)$, respectively. Both are bounded by $a_xa_y+\tau_{x,y}\le C\bigl(f(x)f(y)+A(x,y)\bigr)$ for a constant $C$, which is integrable by \eqref{eq:f-A-integrals}. Dominated convergence gives $I_\wedge\to(\int\varphi)^2>0$ and $I_\triangle\to\int\!\!\int\varphi(x)\varphi(y)\varphi(x-y)$, proving the claim for global clustering.

For the local coefficient, let $U_0$ be the groups joined by the origin, a Poisson process with intensity $\nu\theta g(u)\dd u$. By the probability generating functional,
\[
\E[q(x)]=1-e^{-a_x}=: \varphi_\mu(x),\qquad \operatorname{Cov}(q(x),q(y))=e^{-a_x-a_y}(e^{\tau_{x,y}}-1)\le\tau_{x,y}.
\]
Indeed, $\E[(1-q(x))(1-q(y))]=e^{-a_x-a_y+\tau_{x,y}}$. Hence $\operatorname{Var}(M)\le\int_{\R^{2d}}\tau_{x,y}\dd x\dd y=\nu\theta^3\lVert g\rVert_1^3\to0$, while $\E[M]=\int\varphi_\mu\to\int\varphi=:m>0$ by dominated convergence. Thus $M\to m$ in $L^2$. Moreover, $\operatorname{Var}(q(x))\le\nu\theta^3\lVert g\rVert_2^2\to0$, so $\E|q(x)-\varphi(x)|\to0$ for every $x$. Since $\E|q(x)-\varphi(x)|\le\varphi_\mu(x)+\varphi(x)\le Cf(x)$, dominated convergence gives $\E\int_{\R^d}|q-\varphi|\to0$.

By Lemma~\ref{lem:Q-integral},
\[
QM^2=\int_{\R^{2d}}\Bigl[(1-s(x,y))e^{-\eta(x,y)}+q(x)q(y)\bigl(1-e^{-\eta(x,y)}\bigr)\Bigr]\dd x\dd y,
\]
where $\E[1-s(x,y)]=1-e^{-\tau_{x,y}}\le\tau_{x,y}$ and $\eta(x,y)=a_{x-y}-\tau_{x,y}$. Since $\tau_{x,y}\le\theta\kappa_2$ and $f\le\lVert g\rVert_2^2$, we have $|\eta(x,y)-\kappa f(x-y)/\lVert g\rVert_2^2|\le|\kappa_2-\kappa|+\theta\kappa_2=: \varepsilon_\mu\to0$ uniformly. Therefore,
\[
\E\left|QM^2-\int_{\R^{2d}}q(x)q(y)\varphi(x-y)\dd x\dd y\right|\le\nu\theta^3\lVert g\rVert_1^3+\varepsilon_\mu\E[M^2]\longrightarrow0.
\]
Furthermore, $|q(x)q(y)-\varphi(x)\varphi(y)|\le q(x)|q(y)-\varphi(y)|+\varphi(y)|q(x)-\varphi(x)|$, so
\[
\left|\int_{\R^{2d}}q(x)q(y)\varphi(x-y)\dd x\dd y-\int_{\R^{2d}}\varphi(x)\varphi(y)\varphi(x-y)\dd x\dd y\right|\le(M+m)\int_{\R^d}|q-\varphi|\longrightarrow0
\]
in probability. Since $M^2\to m^2>0$ in probability, it follows that $Q\to C_{\mathrm{RCM}}(\varphi)$ in probability. Corollary~\ref{cor:concentrated}, with normalization $b_\mu=1$, now gives the limit of average local clustering for every choice of the rescaled vertex intensities $\lambda r^d$.

\emph{(c) Saturation.} Since $\nu\theta^2=\kappa_2/\lVert g\rVert_2^2\to\infty$ and $f(z)>0$ exactly when $z\in B$, we have $a_z\to\infty$ for every $z\in B$, whereas $a_z=0$ outside $B$. The wedge probability vanishes unless $x,y\in B$, and on $B\times B$ the union bound gives
\[
1-p_\wedge(x,0,y)\le e^{-a_x}+e^{-a_y}\longrightarrow0.
\]
Likewise, the triangle probability vanishes unless $x,y,x-y\in B$, and whenever these three displacements belong to $B$,
\[
1-p_\triangle(x,0,y)\le e^{-a_x}+e^{-a_y}+e^{-a_{x-y}}\longrightarrow0.
\]
Both probabilities are bounded by the corresponding geometric indicators, so dominated convergence yields $I_\wedge\to|B|^2$ and $I_\triangle\to|B|^2C_{\mathrm{RGG}}$. This proves the global clustering limit.

For the local coefficient, $0\le q\le\1_B$ and $\E\int_B(1-q(x))\dd x=\int_Be^{-a_x}\dd x\to0$. Thus $M\to|B|$ in $L^1$, hence also in $L^2$ since $M\le|B|$. Moreover, $QM^2$ is the integral of the conditional triangle probability given $U_0$, and therefore
\[
0\le QM^2\le\int_{B\times B}\1_{\{\lVert x-y\rVert<2R\}}\dd x\dd y=|B|^2C_{\mathrm{RGG}}.
\]
By the conditional representation of Section~\ref{sec:clustering}, $\E[QM^2]=I_\triangle\to|B|^2C_{\mathrm{RGG}}$. Hence $QM^2\to|B|^2C_{\mathrm{RGG}}$ in $L^1$, and consequently $Q\to C_{\mathrm{RGG}}$ in probability. Applying Corollary~\ref{cor:concentrated} with $b_\mu=1$ completes the proof for arbitrary vertex intensities.
\end{proof}

For bounded-support kernels, Theorem~\ref{thm:boundary} describes global clustering in all three asymptotic boundary regimes of $\mathcal D_g$. Case (a) also covers regimes outside the benchmark scalings. If $\kappa_2\to0$ and $\kappa_1\to\ell\in[0,\infty)$, then $C_{\mathrm{glob},\mu}\to1/(1+\ell)$, independently of the reference kernel and the dimension. In this regime, closed wedges are asymptotically dominated by clique triangles: two neighbours of a vertex are connected essentially when they share a group with that vertex. If instead $\kappa_1\to\infty$ and $\kappa_2\to0$, then $C_{\mathrm{glob},\mu}\to0$. In this latter regime, the bounds in (a) do not by themselves determine the relative proportions of clique and cycle triangles.

\subsection{Clustering in the benchmark limiting regimes}\label{sec:classical-clustering}

The benchmark scalings of Table~\ref{tab:regimes} are paths in the domain $\mathcal D_g$, so their clustering limits follow from Proposition~\ref{prop:two-parameter} and Theorem~\ref{thm:boundary}.

\begin{corollary}[Benchmark regimes]\label{cor:classical-clustering}
Under \eqref{eq:g-assumption}, with fixed vertex intensity $\lambda>0$, consider the scalings of Table~\ref{tab:regimes}.
\begin{enumerate}
\item[(i)] \emph{RCM:} $C_{\mathrm{glob},\mu}$ and $C_{\mathrm{loc},\mu}(\lambda)$ converge to $C_{\mathrm{RCM}}(\varphi)$, with $\varphi=1-e^{-f}$.
\item[(ii)] \emph{RGG:} if $g$ has bounded support, both coefficients converge to $C_{\mathrm{RGG}}$, for every $a\in[0,1/2)$.
\item[(iii)] \emph{RIG:} for every $\mu>0$, $C_{\mathrm{glob},\mu}=\E[Q_*M_*^2]/\E[M_*^2]$, where $(M_*,Q_*)$ are the conditional quantities of the GRIG with group intensity $\alpha$ and kernel $g$; moreover, $C_{\mathrm{loc},\mu}(\lambda)\to\E[Q_*\mid M_*>0]$ as $\mu\downarrow0$.
\item[(iv)] \emph{ER:} both coefficients converge to $C_{\mathrm{RCM}}(\varphi)$, with $\varphi=1-e^{-f}$, and $C_{\mathrm{RCM}}(\varphi)<p_{\mathrm{ER}}=1-e^{-\lVert g\rVert_2^2}$.
\item[(v)] \emph{Complete graph:} if $g$ has bounded support,
both coefficients converge to $C_{\mathrm{RGG}}$, although the
spatially local limit is the complete graph.
\end{enumerate}
\end{corollary}

\begin{proof}
In cases (i) and (iv), Table~\ref{tab:regimes} gives $\kappa_1\to\infty$ and $\kappa_2=\lVert g\rVert_2^2$, so Theorem~\ref{thm:boundary}(b) applies with $\kappa=\lVert g\rVert_2^2$. In case (ii), $\kappa_2=\mu^{1-2a}\lVert g\rVert_2^2\to\infty$, and Theorem~\ref{thm:boundary}(c) applies. In case (iii), $\kappa_1$ and $\kappa_2$ do not depend on $\mu$, so $C_{\mathrm{glob},\mu}$ equals the global clustering coefficient of the GRIG with group intensity $\alpha$ and kernel $g$, which is given by \eqref{eq:global-group-representation}. By Corollary~\ref{cor:clustering-scaling}, $C_{\mathrm{loc},\mu}(\lambda)$ equals the average local clustering coefficient of the same GRIG at vertex intensity $\lambda r^d=\lambda\alpha/\mu\to\infty$, and Corollary~\ref{cor:local-limits} gives its limit.

To prove the inequality in (iv), note that by the Cauchy--Schwarz inequality, $f(x)\le\lVert g\rVert_2^2=f(0)$, and that inequality is strict for $x\ne0$, since equality would force $g$ to agree almost everywhere with a nonzero translate of itself, which is impossible for a nonzero integrable function. 
Hence $\varphi(x)<\varphi(0)=p_{\mathrm{ER}}$ for $x\ne0$, and
\[
\int_{\R^{2d}}\varphi(x)\varphi(y)\varphi(x-y)\dd x\dd y<p_{\mathrm{ER}}\left(\int_{\R^d}\varphi\right)^2,
\]
because the diagonal $\{x=y\}$ is a null set and $\varphi(x)\varphi(y)>0$ on a set of positive measure.

For case~(v), the scaling in Table~\ref{tab:regimes} gives
$\kappa_2=\mu^2\lVert g\rVert_2^2\to\infty$.
If $g$ has bounded support, Theorem~\ref{thm:boundary}(c),
including its conclusion for local clustering, therefore yields
\[
C_{\mathrm{glob},\mu}\longrightarrow C_{\mathrm{RGG}},
\qquad
C_{\mathrm{loc},\mu}(\lambda)\longrightarrow C_{\mathrm{RGG}}.
\]
\end{proof}

In cases (i) and (ii), the limits are the clustering coefficients
of the limiting graphs of Theorem~\ref{thm:local-limits}.
In contrast, cases (iii)--(v) show that spatially local convergence
need not preserve clustering.
In the ER regime, the gap is strict.
In the RIG regime, take $g=\1_{B_R(0)}$: in the limiting RIG,
every group in $B_R(0)$ is joined by all vertices and all other
groups by none, so on any finite vertex set the limiting graph
is either complete or empty, with clustering equal to one
whenever defined.
In the GRIG, instead, wedges whose endpoints lie beyond the
maximum connection distance occur with positive probability
and cannot be closed, so $C_{\mathrm{glob},\mu}<1$.
In case~(v), for bounded-support kernels, both clustering
coefficients converge to $C_{\mathrm{RGG}}<1$, whereas finite
restrictions of the limiting complete graph have clustering
equal to one whenever defined.
Here the geometry persists at the expanding neighbourhood scale,
although every fixed observation window becomes complete.

\begin{table}[htb]
\centering
\scriptsize
\renewcommand{\arraystretch}{1.3}
\begin{tabular}{@{}llccc@{}}
\toprule
& &
$\theta=1$, $\mu r^d=\alpha$ &
$\theta\to0$, $\kappa_2\to\kappa\in(0,\infty)$ &
$\kappa_2\to\infty$\\
\midrule
$r=1$
& spatially local limit
& GRIG
& RCM, $\varphi$
& RGG\\
& global clustering
& $C_{\mathrm{glob},*}$
& $C_{\mathrm{RCM}}(\varphi)$ \ $\checkmark$
& $C_{\mathrm{RGG}}$ \ $\checkmark$\\
& average local clustering
& $C_{\mathrm{loc},*}(\lambda)$
& $C_{\mathrm{RCM}}(\varphi)$ \ $\checkmark$
& $C_{\mathrm{RGG}}$ \ $\checkmark$\\
\midrule
$r\to\infty$
& spatially local limit
& RIG
& ER, $p$
& complete graph\\
& global clustering
& $\displaystyle
   \frac{\E[Q_*M_*^2]}{\E[M_*^2]}$ \ $\times$
& $C_{\mathrm{RCM}}(\varphi)<p$ \ $\times$
& $C_{\mathrm{RGG}}<1$ \ $\times$\\
& average local clustering
& $\E[Q_*\mid M_*>0]$ \ $\times$
& $C_{\mathrm{RCM}}(\varphi)<p$ \ $\times$
& $C_{\mathrm{RGG}}<1$ \ $\times$\\
\bottomrule
\end{tabular}
\caption{
Spatially local limits and corresponding clustering values or
limits for $h_{\theta,r}=\theta g(\cdot/r)$, with fixed
$\lambda>0$ and $\kappa_2=\mu\lVert h_{\theta,r}\rVert_2^2$
(Theorems~\ref{thm:local-limits} and~\ref{thm:boundary},
Corollary~\ref{cor:classical-clustering}).
A checkmark denotes agreement with the clustering coefficient
of the spatially local limit; a cross indicates that agreement
fails in general.
In the first column, $\alpha>0$ is fixed. The starred quantities
refer to the GRIG with group intensity $\alpha$ and kernel $g$;
in the lower block, $r=(\alpha/\mu)^{1/d}$ and $\mu\downarrow0$.
Its global clustering expression holds for every $\mu>0$.
In the other columns, $\mu\to\infty$,
$\varphi=1-e^{-\kappa f/\lVert g\rVert_2^2}$ and
$p=1-e^{-\kappa}$. The last column assumes bounded support of
$g$, except for the spatially local complete-graph limit.
}
\label{tab:regimes-summary}
\end{table}

The distinction between the two rows of \cref{tab:regimes-summary} can be understood through the spatial rescaling of Lemma~4.3. Under this rescaling, a fixed observation window $W$ becomes $W/r$, whereas the clustering coefficients still depend on the entire neighbourhood of a typical vertex in the rescaled graph. When $r\to\infty$, the rescaled window shrinks to a point. The geometry of the full neighbourhood can therefore continue to influence clustering even when the limiting graph on every fixed window is insensitive to vertex positions.

\subsection{The clique--cycle transition in the RGG regime}
The preceding comparison concerns the values of the clustering
coefficients. We now examine a finer question: how triangles are
formed when those coefficients converge to $C_{\mathrm{RGG}}$.
Although the total clustering has the same limit throughout the
RGG scaling $a\in[0,1/2)$, the relative contributions of clique
and cycle triangles undergo a transition at $a=1/3$.
Let $g$ have bounded support, let $B=B_{2R}(0)$, and set
\[
\hat A(x,y):=\frac{A(x,y)}{\lVert g\rVert_3^3}\in[0,1],
\]
where $\hat A\le1$ by H\"older's inequality. For kernels $h_{\theta,r}=\theta g$, the mean number of groups joining vertices at $0$, $x$ and $y$ is $\mu\theta^3A(x,y)=\kappa_3\hat A(x,y)$. For $\kappa\in[0,\infty]$ define
\[
\sigma(\kappa):=\frac{1}{|B|^2}\int_{B\times B}\left(1-e^{-\kappa\hat A(x,y)}\right)\dd x\dd y,
\qquad\text{with }1-e^{-\infty\cdot\hat A}:=\1_{\{A>0\}},
\]
and $\sigma_d:=\sigma(\infty)$. Since $g$ is radial and non-increasing, $\sigma_d$ depends only on the dimension: it is the probability that three balls of radius $R$ centred at $0$, $X$ and $Y$, with $X,Y$ independent and uniform in $B$, have a common intersection of positive volume.

\begin{theorem}[Clique--cycle transition]\label{thm:rgg-clique-cycle}
Under \eqref{eq:g-assumption}, let $g$ have bounded support and consider kernels $h_{\theta,1}=\theta g$, $\theta\in(0,1]$, with fixed vertex intensity $\lambda>0$, such that $\kappa_2\to\infty$ and $\kappa_3\to\kappa_\infty\in[0,\infty]$. Then $C_{\mathrm{glob},\mu}$ and $C_{\mathrm{loc},\mu}(\lambda)$ converge to $C_{\mathrm{RGG}}$, while
\[
C_{\mathrm{glob},\mu}^{(\mathrm{cl})}\longrightarrow\sigma(\kappa_\infty),\qquad C_{\mathrm{loc},\mu}^{(\mathrm{cl})}(\lambda)\longrightarrow\sigma(\kappa_\infty).
\]
The corresponding cycle contributions converge to $C_{\mathrm{RGG}}-\sigma(\kappa_\infty)$. In particular, under the RGG scaling $\theta=\mu^{-a}$, $a\in[0,1/2)$,
\begin{equation}\label{eq:rgg-clique-regimes}
\lim_{\mu\to\infty}C_{\mathrm{glob},\mu}^{(\mathrm{cl})}
=
\begin{cases}
\displaystyle \sigma_d, & 0\le a<1/3,\\[2mm]
\displaystyle \frac{1}{|B|^2}\int_{B\times B}(1-e^{-A(x,y)})\dd x\dd y, & a=1/3,\\[3mm]
\displaystyle 0, & 1/3<a<1/2.
\end{cases}
\end{equation}
\end{theorem}

\begin{proof}
The limits of the total coefficients are given by Theorem~\ref{thm:boundary}(c). By Proposition~\ref{prop:local}, the probability that the vertices at $0,x,y$ form a clique triangle is $p_{\mathrm{cl}}=1-e^{-\kappa_3\hat A(x,y)}$, which converges pointwise to $1-e^{-\kappa_\infty\hat A(x,y)}$. Since $A(x,y)>0$ implies $x,y\in B$, these probabilities vanish outside $B\times B$ and are bounded by one, so dominated convergence gives $I_{\mathrm{cl}}\to|B|^2\sigma(\kappa_\infty)$. The proof of Theorem~\ref{thm:boundary}(c) shows that $I_\wedge\to|B|^2$, which gives the limit of $C_{\mathrm{glob},\mu}^{(\mathrm{cl})}$; the cycle part follows from $I_{\mathrm{cyc}}=I_\triangle-I_{\mathrm{cl}}$.

For the local coefficients, the same proof shows that $M\to|B|$ in $L^2$, with $0\le M\le|B|$. By Lemma~\ref{lem:Q-integral} and the definition of $Q_{\mathrm{cl}}$, $Q_{\mathrm{cl}}M^2=\int\!\!\int(1-s(x,y))\dd x\dd y$, and $\E[1-s(x,y)]$ is the probability that $x$ and $y$ share a group joined by the origin, that is, $p_{\mathrm{cl}}$. Hence $\E[Q_{\mathrm{cl}}M^2]=I_{\mathrm{cl}}$. Since $0\le Q_{\mathrm{cl}}\le1$,
\[
\bigl|\,|B|^2\E[Q_{\mathrm{cl}}]-I_{\mathrm{cl}}\bigr|\le\E\bigl|M^2-|B|^2\bigr|\longrightarrow0,
\]
so $\E[Q_{\mathrm{cl}}]\to\sigma(\kappa_\infty)$. Finally, $w$ is Lipschitz, so $w(\lambda M)\to w(\lambda|B|)>0$ in $L^1$, and
\[
C_{\mathrm{loc},\mu}^{(\mathrm{cl})}(\lambda)=\frac{\E[Q_{\mathrm{cl}}w(\lambda M)]}{\E[w(\lambda M)]}
\]
has the same limit as $\E[Q_{\mathrm{cl}}]$. The cycle part follows from $C_{\mathrm{loc},\mu}(\lambda)=C_{\mathrm{loc},\mu}^{(\mathrm{cl})}(\lambda)+C_{\mathrm{loc},\mu}^{(\mathrm{cyc})}(\lambda)$.

Under the RGG scaling, $\kappa_2=\mu^{1-2a}\lVert g\rVert_2^2\to\infty$ and $\kappa_3=\mu^{1-3a}\lVert g\rVert_3^3$, which tends to $\infty$, $\lVert g\rVert_3^3$ or $0$ according as $a<1/3$, $a=1/3$ or $a>1/3$; for $a=1/3$, $\kappa_\infty\hat A=A$.
\end{proof}

The theorem reveals a transition in the mechanism that produces triangles, one that is invisible at the level of the projected graph. When $\kappa_3\to\infty$, common groups close every triangle for which a common group location is geometrically possible. When $\kappa_3\to0$, such groups disappear, and asymptotically all closed wedges are cycle triangles. In the critical case, both mechanisms contribute, and the limiting decomposition retains information about the shape of $g$ inside its support. Thus, even convergence to the same geometric graph with the same clustering coefficients does not determine the mechanism of triangle formation.

\begin{remark}[The constant $\sigma_d$]\label{rem:sigma}
Since $A(x,y)>0$ requires all pairwise distances among $0,x,y$ to be less than $2R$, we have $\sigma_d\le C_{\mathrm{RGG}}$. In dimension one, three intervals that intersect pairwise have a common point, so $\sigma_1=C_{\mathrm{RGG}}=3/4$. For $d\ge2$, the inequality is strict: if $0,x,y$ form an equilateral triangle with side length in $(\sqrt3R,2R)$, the three balls intersect pairwise, but their common intersection is empty, because the circumradius of the triangle exceeds $R$; both conditions are open, so they hold on a set of positive measure.
\end{remark}

\begin{remark}[Correlations and higher-order groups]\label{rem:higher-order}
By Proposition~\ref{prop:edge-correlations}, for kernels $\theta g$ the non-edge events of a wedge satisfy
\[
\frac{\Prob(0\nsim x,\ 0\nsim y)}{\Prob(0\nsim x)\,\Prob(0\nsim y)}=e^{\kappa_3\hat A(x,y)}.
\]
Under the RGG scaling, this ratio tends to $\infty$, $e^{A(x,y)}$ or $1$ according as $a<1/3$, $a=1/3$ or $a>1/3$, for every $(x,y)$ with $A(x,y)>0$. Thus the transition of Theorem~\ref{thm:rgg-clique-cycle} is also visible in the dependence between non-edges at finite $\mu$. More generally, groups joined by $s$ close vertices have mean approaching $\kappa_s=\mu^{1-sa}\lVert g\rVert_s^s$ as their positions coalesce under the RGG scaling, which suggests analogous thresholds $a=1/s$ for larger cliques; we do not pursue this here.
\end{remark}

\section{Conclusion and open directions}\label{sec:conclusion}

In this work, we investigated the subgraph structure and clustering properties of geometric random intersection graphs, with particular emphasis on the interplay between spatial geometry and group-mediated dependencies. Starting from an exact representation of subgraph probabilities through independent group types, we derived expressions for global and average local clustering and characterized their dependence on the vertex intensity. We then established spatially local convergence to several classical random graph models, identifying distinct parameter scalings leading to random intersection graphs, random connection models, random geometric graphs, Erd\H{o}s--R\'enyi graphs, and the complete graph.

A central conclusion is that the convergence of random graph models on bounded spatial windows does not necessarily imply the convergence of their clustering coefficients. While clustering is preserved in the RCM and RGG regimes, the RIG, ER, and complete-graph limits demonstrate that contributions from increasingly distant vertices may remain relevant even when the graph converges locally. In particular, the ER and complete-graph regimes provide strict examples in which the parameter limit and the infinite-volume limit do not commute. Our general clustering convergence criterion provides a unified framework for understanding these different behaviours.

Beyond the limiting values of clustering, our results also reveal that the same limiting graph can arise from substantially different mechanisms of triangle formation. In the RGG regime, the total clustering converges to a geometric constant depending only on the dimension, while the relative contributions of clique and cycle triangles exhibit a transition at the scaling exponent $a=1/3$. Thus, neither the limiting graph nor its clustering coefficient alone fully captures the structure of the underlying group interactions.

Several questions remain open. One natural direction is to understand the role of the tail behaviour of the membership kernel in regimes where the spatial range of connections diverges. In particular, it would be interesting to determine whether rapidly decaying and regularly varying kernels lead to distinct universality classes for clustering. Another direction concerns joint scalings of the vertex and group intensities, which may reveal additional asymptotic regimes and further instances of non-commuting limits. Finally, the distinction between spatially local convergence and the convergence of infinite-volume observables suggests investigating analogous questions for more general subgraph densities and higher-order clustering statistics.


\subsection*{Acknowledgments}
A. Bianchi was partially supported by the GNAMPA project ``Ising model on random graphs with community structure: equilibrium and dynamical transitions''. R. Michielan was partially supported by the GAMING project ``Graph Algorithms and MinINg for Green agents'' within the FAIR partnership ``Future Artificial Intelligence Research'', funded by the European Union through NextGenerationEU under the Italian National Recovery and Resilience Plan (PNRR).

\subsection*{Data availability}
No datasets were generated or analysed in this study.

\subsection*{Conflict of interest}
The authors declare no conflicts of interest.

\bibliographystyle{abbrv}
\bibliography{references}

@article{deijfen-michielan,
  author  = {Deijfen, Maria and Michielan, Riccardo},
  title   = {Geometric random intersection graphs with general connection probabilities},
  journal = {Journal of Applied Probability},
  year    = {2024},
  volume  = {61},
  number  = {4},
  pages   = {1343--1360},
  doi     = {10.1017/jpr.2024.18}
}

@article{bloznelis,
  author  = {Bloznelis, Mindaugas},
  title   = {Degree and clustering coefficient in sparse random intersection graphs},
  journal = {The Annals of Applied Probability},
  year    = {2013},
  volume  = {23},
  number  = {3},
  pages   = {1254--1289},
  doi     = {10.1214/12-AAP874}
}

@book{last-penrose,
  author    = {Last, G{\"u}nter and Penrose, Mathew},
  title     = {Lectures on the {Poisson} Process},
  publisher = {Cambridge University Press},
  year      = {2017},
  doi       = {10.1017/9781316104477},
  isbn      = {9781107088016}
}

@book{penrose,
  author    = {Penrose, Mathew},
  title     = {Random Geometric Graphs},
  publisher = {Oxford University Press},
  year      = {2003},
  doi       = {10.1093/acprof:oso/9780198506263.001.0001},
  isbn      = {9780198506263}
}

@article{karonski-scheinerman-singer,
  author  = {Karo{\'n}ski, Micha{\l} and Scheinerman, Edward R. and Singer-Cohen, Karen B.},
  title   = {On random intersection graphs: the subgraph problem},
  journal = {Combinatorics, Probability and Computing},
  year    = {1999},
  volume  = {8},
  number  = {1-2},
  pages   = {131--159},
  doi     = {10.1017/S0963548398003459}
}

@article{rybarczyk-stark,
  author  = {Rybarczyk, Katarzyna and Stark, Dudley},
  title   = {{Poisson} approximation of counts of induced subgraphs in random intersection graphs},
  journal = {Discrete Mathematics},
  year    = {2017},
  volume  = {340},
  number  = {9},
  pages   = {2183--2193},
  doi     = {10.1016/j.disc.2017.04.013}
}

@article{bachmann-reitzner,
  author  = {Bachmann, Sascha and Reitzner, Matthias},
  title   = {Concentration for {Poisson} {$U$}-statistics: subgraph counts in random geometric graphs},
  journal = {Stochastic Processes and their Applications},
  year    = {2018},
  volume  = {128},
  number  = {10},
  pages   = {3327--3352},
  doi     = {10.1016/j.spa.2017.11.001}
}

@article{hofstad-komjathy-vadon,
  author  = {van der Hofstad, Remco and Komj{\'a}thy, J{\'u}lia and Vadon, Vikt{\'o}ria},
  title   = {Random intersection graphs with communities},
  journal = {Advances in Applied Probability},
  year    = {2021},
  volume  = {53},
  number  = {4},
  pages   = {1061--1089},
  doi     = {10.1017/apr.2021.12}
}

@article{hofstad-hoorn-maitra,
  author  = {van der Hofstad, Remco and van der Hoorn, Pim and Maitra, Neeladri},
  title   = {Scaling of the Clustering Function in Spatial Inhomogeneous Random Graphs},
  journal = {Journal of Statistical Physics},
  year    = {2023},
  volume  = {190},
  number  = {6},
  pages   = {110},
  doi     = {10.1007/s10955-023-03122-6}
}

@mastersthesis{de-nicolo-thesis,
  author = {De Nicol{\`o}, Alberto},
  title  = {The geometric random intersection graph model: clustering and local Structure},
  school = {Universit{\`a} degli Studi di Padova},
  year   = {2026},
  note   = {Academic year 2025/2026},
  url    = {https://hdl.handle.net/20.500.12608/115675}
}

@article{rybarczyk-serafin-subgraph-probabilities,
  title={Subgraph probabilities in random intersection graphs},
  author={Rybarczyk, Katarzyna and Serafin, Grzegorz},
  journal={arXiv preprint arXiv:2610.27117},
  year={2026}
}

@article{NguyenZessin79,
  author  = {Nguyen, Xuan Xanh and Zessin, Hans},
  title   = {Ergodic theorems for spatial processes},
  journal = {Zeitschrift f{\"u}r Wahrscheinlichkeitstheorie und Verwandte Gebiete},
  volume  = {48},
  number  = {2},
  pages   = {133--158},
  year    = {1979},
  doi     = {10.1007/BF01886869}
}

@article{last-nestmann-schulte,
  author        = {Last, G{\"u}nter and Nestmann, Franz and Schulte, Matthias},
  title         = {The random connection model and functions of edge-marked {Poisson} processes: second order properties and normal approximation},
  journal       = {The Annals of Applied Probability},
  year          = {2021},
  volume        = {31},
  number        = {1},
  pages         = {128--168},
  doi           = {10.1214/20-AAP1585},
  eprint        = {1808.01203},
  archivePrefix = {arXiv},
  url           = {https://arxiv.org/abs/1808.01203}
}

@article{benjamini-schramm,
  author  = {Benjamini, Itai and Schramm, Oded},
  title   = {Recurrence of Distributional Limits of Finite Planar Graphs},
  journal = {Electronic Journal of Probability},
  year    = {2001},
  volume  = {6},
  number  = {23},
  pages   = {1--13},
  doi     = {10.1214/EJP.v6-96}
}

@article {brun26,
    AUTHOR = {Brun, Morten and Hirsch, Christian and Juh\'asz, P\'eter and
              Otto, Moritz},
     TITLE = {Random connection hypergraphs},
   JOURNAL = {Stochastic Models},
  FJOURNAL = {Stochastic Models},
    VOLUME = {42},
      YEAR = {2026},
    NUMBER = {3},
     PAGES = {345--387},
      ISSN = {1532-6349,1532-4214},
   MRCLASS = {05C80 (05C65 60C05 60F05)},
  MRNUMBER = {5106509},
       DOI = {10.1080/15326349.2026.2651174},
       URL = {https://doi.org/10.1080/15326349.2026.2651174},
}

@article{lichev2023,
  author  = {Lichev, Lyuben and Lodewijks, Bas and Mitsche, Dieter and Schapira, Bruno},
  title   = {On the first and second largest components in the percolated random geometric graph},
  journal = {Stochastic Processes and their Applications},
  volume  = {164},
  pages   = {311--336},
  year    = {2023},
  doi     = {10.1016/j.spa.2023.07.008}
}

@book{meester-roy,
  author = {Meester, Ronald and Roy, Rahul},
  title = {Continuum Percolation},
  publisher = {Cambridge University Press},
  year = {1996}
}

@article{deijfen-kets,
  author = {Deijfen, Maria and Kets, Willemien},
  title = {Random intersection graphs with tunable degree distribution and
    clustering},
  journal = {Probability in the Engineering and Informational Sciences},
  volume = {23},
  number = {4},
  pages = {661--674},
  year = {2009}
}

@article{fill-scheinerman-singer,
  author = {Fill, James Allen and Scheinerman, Edward R. and
    Singer-Cohen, Karen B.},
  title = {Random intersection graphs when $m=\omega(n)$: an equivalence
    theorem relating the evolution of the {$G(n,m,p)$} and {$G(n,p)$}
    models},
  journal = {Random Structures \& Algorithms},
  volume = {16},
  number = {2},
  pages = {156--176},
  year = {2000}
}

@article{rybarczyk-equivalence,
  author = {Rybarczyk, Katarzyna},
  title = {Equivalence of a random intersection graph and {$G(n,p)$}},
  journal = {Random Structures \& Algorithms},
  volume = {38},
  number = {1--2},
  pages = {205--234},
  year = {2011},
  doi = {10.1002/rsa.20356}
}

@article{kim-lee-na,
  author  = {Kim, Jeong Han and Lee, Sang June and Na, Joohan},
  title   = {On the total variation distance between the binomial random graph and the random intersection graph},
  journal = {Random Structures \& Algorithms},
  year    = {2018},
  volume  = {52},
  number  = {4},
  pages   = {662--679},
  doi     = {10.1002/rsa.20750}
}

@article{brennan-bresler-nagaraj,
  author  = {Brennan, Matthew and Bresler, Guy and Nagaraj, Dheeraj},
  title   = {Phase transitions for detecting latent geometry in random graphs},
  journal = {Probability Theory and Related Fields},
  year    = {2020},
  volume  = {178},
  number  = {3--4},
  pages   = {1215--1289},
  doi     = {10.1007/s00440-020-00998-3}
}

\end{document}